\input amstex
\documentstyle{amsppt}

\input xy  
\xyoption{all}


\define\mathbb{\Bbb}
\define\mathscr{\Cal}
\define\mathcal{\Cal}
\define\mathfrak{\frak}
\define\em{\it}

\newdimen\jedn
\catcode`@=11
\newcount\@multicnt
\newdimen\@xdim
\newdimen\@ydim
 
\long\def\put(#1,#2)#3{\raise#2\jedn\hbox
to\z@{\kern#1\jedn#3\hss}\ignorespaces
}
\def\grub#1 {\special{em:linewidth #1}}
\def\MT #1 #2 {\put(#1,#2){\special{em:moveto}}}
\def\LT #1 #2 {\put(#1,#2){\special{em:lineto}}}
 
\long\def\multi(#1,#2)(#3,#4)#5#6{\@multicnt=#5\relax
\@xdim=#1\jedn \@ydim=#2\jedn
\loop\ifnum \@multicnt>\z@
\raise\@ydim\hbox
to\z@{\kern\@xdim\hbox{#6}\hss}\advance\@multicnt by-1 \advanc
e\@xdim
#3\jedn\advance\@ydim #4\jedn \repeat \ignorespaces}
 
\grub0.4pt
\jedn=0.8mm
\catcode`@=12
 
\long\def\cput(#1,#2,#3){\put(#1,#2){\hbox to0pt{\hss$#3$\hss}}}
 
\long\def\rput(#1,#2,#3){\put(#1,#2){\hbox to0pt{$#3$\hss}}}
\long\def\lput(#1,#2,#3){\put(#1,#2){\hbox to0pt{\hss$#3$}}}

\NoBlackBoxes

\topmatter
\title
Eigenvalues, invariant factors, highest weights, and Schubert calculus 
\endtitle
\rightheadtext{Eigenvalues and Schubert calculus}
\author William Fulton \endauthor
\address  University of Michigan,  
	Ann Arbor, MI 48109-1109  \endaddress
\email wfulton\@math.lsa.umich.edu \endemail
\date March 25, 2000 \enddate
\thanks The author was partly supported by NSF Grant \#DMS9970435 
\endthanks

\abstract
We describe recent work of Klyachko, Totaro, Knutson, and Tao, that 
characterizes eigenvalues of sums of Hermitian matrices, and decomposition 
of tensor products of representations of $GL_n(\mathbb{C})$.  We explain 
related applications to invariant factors of products of matrices,  
intersections in Grassmann varieties, and singular values of sums and 
products of arbitrary matrices.   
\endabstract

\toc
\widestnumber\head{10}	
\head   1.  Eigenvalues of sums of Hermitian and real symmetric matrices 
\endhead
\head   2.  Invariant factors \endhead
 \head  3.  Highest weights \endhead
 \head  4.  Schubert calculus \endhead
 \head  5.  Singular values of sums and products \endhead
 \head  6.  First steps toward the proofs \endhead
 \head  7.  Filtered vector spaces, geometric invariant theory, and stability \endhead
\head 8. Coadjoint orbits and moment maps \endhead
 \head  9.  Saturation \endhead
 \head  10.  Proofs of the theorems \endhead
 \head  11. Final remarks \endhead
\endtoc

\endtopmatter
\document

Recent breakthroughs, primarily by A. Klyachko, B. Totaro, A. Knutson, and 
T. Tao, with contributions by P. Belkale and C. Woodward, have led to 
complete solutions of several old problems involving the various notions in 
the title.  Our aim here is to describe this work, and especially to show how 
these solutions are derived from it.  Along the way, we will see that these 
problems are also related to other areas of mathematics, including 
geometric invariant theory, symplectic geometry, and combinatorics.  
In addition, we present some 
related applications to singular values of arbitrary matrices. 

Although many of the theorems we state here have not appeared elsewhere, 
their proofs are mostly ``soft'' algebra based on the hard geometric or 
combinatorial work of others. Indeed, this paper emphasizes concrete 
elementary arguments.  We do give some new examples and 
counterexamples, and raise some new open questions. 

We have attempted to point to the sources and to some of the key partial 
results that had been conjectured or proved before.  However, there is a very 
large literature, particularly for linear algebra problems about eigenvalues, 
singular values, and invariant factors.  We have listed only a few of these 
articles, from whose bibliographies, we hope, an interested reader can trace 
the history; we apologize to the many whose work is not cited directly. 

We begin in the first five sections by describing each of the problems, 
together with some of their early histories, and we state as theorems the 
new solutions to these problems.  In Section 6 we 
 describe the steps toward these 
solutions that were carried out before the recent breakthroughs.  Then we 
discuss the recent solutions, and explain how these theorems follow from 
the work of the above mathematicians.  Sections 7, 8, 9, and 10 also contain 
variations and generalizations of some of the theorems stated in the first five 
sections, as well as attributions of the theorems to their authors.   

One of our fascinations with this subject, even now that we have proofs of the 
theorems, is the challenge to understand in a deeper way why all these 
subjects are so closely related.  It is a particular challenge in each of 
these areas to understand why the solutions can be described inductively.  

I am grateful to many people for advice and help in preparing this article, 
especially: P. Belkale, A. Buch, L. Chen, 
J. Day, P. Diaconis, A. Knutson, R. Lazarsfeld, C. K. Li,
A. Okounkov, Y. T. Poon,  P. Pragacz, J. F. Queir\'o, P. Sarnak, 
F. Sottile, R. Steinberg, J. Stembridge, T. Y. Tam, 
H. Tamvakis, T. Tao, B. Totaro, C. Woodward, 
A. Zelevinsky, and a referee. 

\head 1.  Eigenvalues of sums of Hermitian and real symmetric matrices 
\endhead

The first problem goes back to the nineteenth century:  What can be said 
about the eigenvalues of a sum of two Hermitian (or real symmetric) 
matrices, in terms of the eigenvalues of the summands?

It is a basic fact of linear algebra that all of the eigenvalues 
of any Hermitian 
or real symmetric matrix are real.  We consider $n$ by $n$  matrices, with 
$n$ fixed.  If $A$ is a real symmetric matrix, its eigenvalues describe the 
quadratic form $q_A(x) = x^t A\, x$ in an appropriate orthogonal 
coordinate system.  For example, if the eigenvalues are positive, the inverses 
of the square roots of the eigenvalues are half the lengths of the principal 
axes of the ellipsoid  $q_A(x) =1$.

We always list the $n$ eigenvalues of such a matrix in decreasing order, 
including any eigenvalue as often as its multiplicity.  If $A$, $B$, and $C$ 
are Hermitian $n$ by $n$ matrices, we denote the eigenvalues of $A$ by 
$$
	\alpha \, :\,\,\,\,
   \alpha_1 \geq  \alpha_2 \geq  \ldots  \geq  \alpha_n,  
$$
and similarly write $\beta$ and $\gamma$  for the eigenvalues 
(spectra) of $B$ and 
$C$.  The question becomes: 
\block
	What $\alpha$, $\beta $, and $\gamma$ can be the eigenvalues of $n$ 
	by $n$ Hermitian (or real symmetric) matrices $A$, $B$, and $C$, 
	with $C = A + B$?   
\endblock

For example, one can fix $A$ and take $B$ small, thus regarding $C$ as a 
deformation of $A$.  In the real symmetric case, one is then asking how the 
shape of the hypersurface $q_C(x) =1$ compares with that of $q_A(x) =1$.

One can, if one wishes, take $A = D(\alpha)$ to be the diagonal matrix with 
entries $\alpha_1, \ldots , \alpha_n$ down the diagonal, and $B = U\, 
D(\beta)\,U^*$, with $U$ a unitary (or orthogonal in the real symmetric 
case) matrix.  We are looking for the eigenvalues of $D(\alpha) + U\, 
D(\beta)\, U^*$ as $U$ varies over the unitary group $U(n)$.  This was the 
approach in much of the classical work on the problem, but it will not 
play an important role  
here.

There is one obvious necessary condition, that the trace of $C$ be the sum 
of the traces of $A$ and $B$:
$$
	\sum_{i=1}^n \gamma_i\,\, =\,\, \sum_{i=1}^n \alpha_i  + \sum_{i=1}^n 
	\beta_i. \tag$*$
$$

There is a long history of results that put additional 
necessary conditions on the 
possible eigenvalues.  The first goes back more than a century, and is a 
reasonable exercise for a linear algebra class:
$$
	\gamma_1 \,\leq\, \alpha_1 + \beta_1.  \tag1
$$

The first significant result was given in 1912 by H. Weyl [W]:
$$
	\gamma_{i+j-1}\, \leq\, \alpha_i + \beta_j  \,\, \text{ whenever } 
\,\, i+j-1 \leq n. \tag2
$$

Here is a typical application of these inequalities. If $A$ and $B$ are 
Hermitian $n$ by $n$ matrices that differ by a matrix of rank at most $r$, 
then their eigenvalues $\alpha$  and $\beta$  satisfy the inequalities 
$$
	\alpha_{k+r} \leq \beta_k \,\text{ and }\, \beta_{k+r} \leq \alpha_k 
\,\,\, \text{ for }\, 1 \leq k, \,\,\,\, k+r \leq n.  
$$
Indeed, these inequalities follow by applying (2) to the triples 
$(B, A-B, A)$ and $(A, B-
A, B)$, with $j = r+1$.  When $r = 1$, this is an {\em interlacing theorem:} 
between any two odd numbered (or even numbered) eigenvalues of $A$ 
there is at least one eigenvalue of $B$.  There are similar results for 
eigenvalues of principal minors of $A$, which we describe in Section 6.

In fact, for $n = 2$, it is not hard to verify directly that the conditions ($*$), 
(1), and (2), which say
$$\align
	&\gamma_1 + \gamma_2 = \alpha_1 + \alpha_2 + \beta_1 + \beta_2 
\\
	&\gamma_1 \leq \alpha_1 + \beta_1, \,\,\,  \gamma_2 \leq 
\alpha_2 + 	\beta_1,  \text { and }\,\,  \gamma_2 \leq \alpha_1 + \beta_2
\endalign$$
are both necessary and sufficient for the existence of $2$ by $2$ Hermitian 
(or real symmetric) matrices with these eigenvalues.

In 1949 K. Fan [F] found some other necessary conditions:
$$
	\sum_{i=1}^r \gamma_i \leq \sum_{i=1}^r \alpha_i  + \sum_{i=1}^r 
\beta_i   \quad \text{ for any } \,\, r < n. \tag3
$$

This question was featured in Gel'fand's seminar in Moscow.  In 1950 V. B. 
Lidskii [L1], cf. [BG], found the following necessary condition.  For this, 
regard an $n$-tuple $\alpha$  of eigenvalues as a point $(\alpha_1, \ldots  , 
\alpha_n)$ in $\mathbb{R}^n$.  This result asserts that the point 
$\gamma$  must be in the convex hull of the points $\alpha  + 
\beta_\sigma$,  where $\sigma$ varies over the symmetric group $S_n$, 
and  $\beta_\sigma$ denotes $(\beta_{\sigma(1)}, \ldots  , 
\beta_{\sigma(n)})$.  Although this looks quite different from (1) -- (3), H. 
Wielandt [Wi] showed that this geometric condition is equivalent to the 
inequalities:
$$
	\sum_{i  \in  I} \gamma_i \leq \sum_{i \in  I} \alpha_i + 
\sum_{i=1}^r \beta_i   \tag4
$$
for every subset $I$ of $\{1, \ldots , n \}$ of cardinality $r$, 
for all $r < n$.  
Of course, the same inequalities are valid when the roles of $\alpha$ and 
$\beta$  are interchanged, and we include them in the list we are making.  
In 1962, A. Horn [H2] showed that the inequalities listed so far, with 
the addition of 
$$
 	\gamma_2 + \gamma_3 \leq \alpha_1 + \alpha_3 + \beta_1 + 
\beta_3, \tag5
$$
are necessary and sufficient for the existence of 
$3$ by $3$ Hermitian matrices with these eigenvalues. 

Other inequalities were found, all having the form
$$
   \sum_{k \in  K} \gamma_k \,\, \leq \,\, \sum_{i \in  I} \alpha_i + \sum_{j \in  J} 
\beta_j,   \tag$*_{IJK}$
$$
for certain subsets $I$, $J$, $K$ of $\{1, \ldots, n\}$ of the same 
cardinality $r$, with $r < n$.  We always write the subsets in increasing 
order, so 
$$I = \{i_1 < i_2 < \ldots < i_r\}, $$
$J = \{j_1 <  \ldots < j_r\}$, and $K = \{k_1 < \ldots < k_r\}$. 

For example, there is a result of L. Freede and R. C. Thompson 
[TF], which generalizes the Weyl, Fan, and Lidskii-Wielandt inequalities:  
$$	
	\sum_{k \in  K} \gamma_k \leq \sum_{i \in  I} \alpha_i + 
	\sum_{j \in  J} \beta_j,\, \text{ if }\, k_p = i_p + j_p - p 
	\,\text{ for }\, 1 \leq p \leq r.   \tag6
$$
One has these inequalities for any $I$ and $J$, provided that $i_r + j_r \leq 
n + r$.

Note that any inequality of the form ($*_{IJK}$) can be subtracted from the 
equality ($*$) to give an inequality 
$$
	\sum_{k \in K^c} \gamma_k \geq  \sum_{i \in I^c} \alpha_i + 
	\sum_{j \in J^c} \beta_j, \tag7
$$
where $I^c$, $J^c$, and $K^c$ are the complements of $I$, $J$, and $K$ 
in $\{1, \ldots ,n\}$.  We do not list any of these. 

In his remarkable paper [H2], Horn undertook a systematic study of such 
inequalities.  In fact, he prescribed sets of triples $(I,J,K)$, and he 
conjectured that the inequalities ($*_{IJK}$) for these triples would give 
both necessary and sufficient conditions for a triple $(\alpha ,\beta ,\gamma 
)$ to arise as eigenvalues of Hermitian matrices $A$, $B$, and $C$ with $C 
= A + B$.  

Horn defined sets $T_r^n$ of triples $(I,J,K)$ of subsets of $\{1, \ldots , n 
\}$ of the same cardinality $r$, by the following inductive procedure.  Set 
$$
	U_r^n = \{(I,J,K) \,\, | \,\, \sum_{i \in I} i + \sum_{j \in J} j = 
	\sum_{k \in K} k + r(r+1)/2\}. \tag8
$$
All the triples that we have listed are in $U_r^n$.  When $r = 1$, set 
$T_1^n = U_1^n$.  (The inequalities specified by $(I,J,K)$ in $T_1^n$ are 
Weyl's inequalities (2).)  In general, 
$$\aligned
	T_r^n = \{(I,J,K) &\in U_r^n \, | \, \text{ for all } p < r  \text{ and 
all } 	(F,G,H) \text{ in } T_p^r,  \\
	 & \sum_{f \in F} i_f +\sum_{g \in G} j_g 	
	\leq \sum_{h \in H} k_h + p(p+1)/2\}. 
\endaligned \tag9 $$
\proclaim{Horn's Conjecture}  A triple $(\alpha,\beta,\gamma)$
occurs as eigenvalues of Hermitian $n$ by $n$ matrices $A$, $B$, $C$  with 
$ C = A+B$ if and only if $\sum \gamma_i = \sum \alpha_i + \sum \beta_i$ 
and the inequalities \rom{(}$*_{IJK}$\rom{)} hold for every $(I,J,K)$ in 
$T_r^n$, for all $r < n$. \endproclaim

B.V. Lidskii, son of V. B. Lidskii, announced a proof of Horn's conjecture in 
1982 [L2], but no details of this have ever appeared.

\proclaim{Theorem 1}  Horn's conjecture is true. \endproclaim

In fact, from the recent work to be discussed below, the surprising 
inductive form of this eigenvalue problem can be given in another way.  To 
state this we use a standard correspondence between finite sets $I = \{i_1 < 
i_2 < \ldots < i_r\}$ of $r$ positive integers, and partitions $\lambda  = 
(\lambda_1 \geq \lambda_2 \geq \ldots  \geq  \lambda_r)$ of nonnegative 
integers of length at most $r$.  This correspondence is obtained by defining
$$
	\lambda  = \lambda(I) = (i_r - r, i_{r-1} - (r-1), \ldots , i_2 - 2, i_1 - 
1). \tag10
$$
We let $J$ correspond to $\mu$, and $K$ correspond to $\nu$, by this 
same recipe.  When we say $(\lambda ,\mu,\nu)$ {\bf corresponds to} 
$(I,J,K)$, we mean $(\lambda ,\mu,\nu) = 
(\lambda(I),\lambda(J),\lambda(K))$.  (To recover the subsets from the 
partitions, the integer $r$ must be specified, but this will always be the case 
for us.)

If $(I,J,K)$ corresponds to $(\lambda ,\mu,\nu)$, note that $(I,J,K)$ is in 
$U_r^n$ exactly when $\sum \nu_i = \sum \lambda_i + \sum \mu_i$.
Note also that each of $\lambda$, $\mu$, and $\nu$ is a sequence of $r$ real 
numbers listed in decreasing order, so it makes sense to ask if {\em they} 
arise as eigenvalues of a triple of Hermitian $r$ by $r$ matrices. 

\proclaim{Theorem 2}  A triple $(I,J,K)$ is in $T_r^n$ if and only if the 
corresponding triple $(\lambda ,\mu,\nu)$ occurs as eigenvalues of a triple 
of Hermitian $r$ by $r$ matrices, with the third the sum of the first two. 
\endproclaim

For example, one can use this theorem 
to find (all) such triples with $r = n-1$.  They 
consist of complements of integers $i$,$j$, and $k$ with $i + j = n + k$.  
The triples $(\lambda ,\mu,\nu)$ corresponding to these subsets consist of 
partitions with only $1$'s and $0$'s, with $n-i$, $n-j$, and $n-k$ $1$'s in 
$\lambda$, $\mu$, and $\nu$ respectively.  These occur as eigenvalues of 
diagonal matrices of size $n-1$ with only $1$'s and $0$'s as entries; the 
first matrix has $n-i$ $1$'s at the beginning of the diagonal, the second 
has $n-
j$ $1$'s at the end, and their sum has a total of $n-k$ $1$'s at the ends.  
From (7) this gives the inequalities
$$
	\gamma_k  \geq  \alpha_i + \beta_j  \quad \text{ if } \quad i + j = n + 
k. \tag11
$$

For another example, the triple $(I,J,K) = (\{1,3,5\}, \{1,3,5\}, \{2,4,6\})$ 
is in $T_3^6$.  In this case the corresponding triple $(\lambda ,\mu,\nu) = 
((2,1,0), (2,1,0), (3,2,1))$ arises from the triple of diagonal $3$ by $3$ 
matrices with diagonal entries $(2,1,0)$, $(1,0,2)$, and $(3,1,2)$. 

We will give several other characterizations of the sets $T_r^n$ in terms of 
Littlewood-Richardson coefficients (see Theorems 11 and 12)
in Section 3.

Horn's paper also contains results and examples about the case of real 
symmetric matrices, but he did not explicitly make his conjecture for this 
case.  In fact, some of the methods that were used to prove inequalities 
broke down in the real case (see Example 1), so there was some doubt 
whether the real symmetric and complex Hermitian cases would coincide; 
a positive answer, however, was conjectured in [DST].  Theorems 1 
and 2 are indeed true in the real symmetric case:

\proclaim{Theorem 3}  A triple $(\alpha ,\beta ,\gamma)$ occurs as 
eigenvalues for a triple of real symmetric matrices
if and only if it occurs for a
triple of complex Hermitian matrices. \endproclaim

There have also been a few cases of inequalites 
($*_{IJK}$) proved for 
Hermitian (self-conjugate) quaternionic matrices (e.g. [X]).  The 
complete answer, as suggested by Steinberg and Totaro, is the same 
as in the real and complex situations: 

\proclaim{Theorem 4}  A triple $(\alpha ,\beta ,\gamma)$ occurs as
eigenvalues for a triple of quaternionic  
Hermitian matrices if and only if it occurs for a
triple of complex Hermitian matrices. \endproclaim

V. B. Lidskii and Horn investigated the situation when some inequality 
($*_{IJK}$) becomes an equality, and they stated and proved special cases of 
the following theorem.

\proclaim{Theorem 5}  If Hermitian matrices $A$, $B$, and $C = A+B$ have 
eigenvalues $\alpha$, $\beta$, and $\gamma$, and some inequality 
\rom{(}$*_{IJK}$\rom{)} occurs with equality for $(I,J,K)$ in $T_r^n$, 
then there is an $r$-dimensional subspace $L$ of $\mathbb{C}^n$ such that 
$A$, $B$, and $C$ map $L$ into itself. \endproclaim

This means that, after a unitary change of coordinates, the three matrices 
have a block diagonal form $\left(\smallmatrix P & 0 \\ 0 & Q 
\endsmallmatrix\right)$ with $P$ $r$ by $r$ and $Q$ $n-r$ by $n-r$ 
Hermitian matrices.

\example{Example 1}  The assertion of Theorem 5 is {\em not} true for 
real symmetric matrices and real subspaces of $\mathbb{R}^n$, although it 
is true for $n \leq 5$.  Here is an explicit example with $n = 6$.  Take $A$ 
to be any diagonal matrix with entries $(x,x,y,y,z,z)$ down the diagonal, with 
$x$, $y$, and $z$ distinct real numbers that sum to zero.  Take $B$ to be 
the matrix
$$
\pmatrix
15 & 0 & -32 & -3 & 35 & -3 \\
0 & 15 & 3 & -32 & 3 & 35 \\
-32 & 3 & -17 & 0 & 51 & 19 \\
-3 & -32 & 0 & -17 & -19 & 51 \\
35 & 3 & 51 & -19 & 2 & 0 \\
-3 & 35 & 19 & 51 & 0 & 2  
\endpmatrix ,
 $$   
which has eigenvalues $(56, 56, 28, 28, -84, -84)$.  Then $C = A+B$ will 
also have eigenvalues consisting of three pairs of distinct numbers that sum 
to zero.  The inequality ($*_{IJK}$), with $I = J = \{1,3,5\}$, and $K = 
\{2,4,6\}$ is an equality, since $\alpha_1 + \alpha_3 + \alpha_5 = 0$, 
$\beta_1 + \beta_3 + \beta_5= 0$, and $\gamma_2 + \gamma_4 + 
\gamma_6 = 0$.  There are exactly two $3$-dimensional subspaces 
preserved by $A$ and $B$ (and therefore $C$), and both are complex; 
letting $e_1, \ldots , e_6$ be the standard basis for $\mathbb{C}^n$, they 
are the span $L$ of $e_1 + ie_2$, $e_3 + ie_5$, and $e_5 + ie_6$, and its 
orthogonal subspace $L^\perp$, spanned by $e_1 - ie_2$, $e_3 - ie_5$, 
and $e_5 - ie_6$.  (The fact that $C$ preserves these subspaces implies that 
its eigenvalues come in pairs as asserted.)  We will see some explanation 
for this example later (cf. [TT2]). \endexample
  
It should be understood that one is looking for a minimal, or at least a small 
set of inequalities ($*_{IJK}$).  Any such inequality determines many others, 
and one usually does not want to list these trivial consequences.  For 
example, from (5) one can deduce immediately that $\gamma_2 + 
\gamma_4 \leq \alpha_1 + \alpha_3 + \beta_1 + \beta_2$, simply because 
$\gamma_4 \leq \gamma_3$ and $\beta_3 \leq \beta_2$. 
As A. Buch points out, there is a special situation for $n = 2$, as the 
inequalities $\alpha_1 \geq \alpha_2$, $\beta_1 \geq \beta_2$, and 
$\gamma_1 \geq  \gamma_2$ follow from the equality ($*$) and  
the three inequalities ($*_{IJK}$).
With this exception, the 
inequalities in Horn's conjecture are in fact minimal for $n \leq 5$, which is 
the region that Horn investigated most thoroughly.    We will see, however, 
that they are {\em not} minimal for larger $n$, although that was thought 
to be the case until quite recently.  

Although the triple $(I,J,K) = (\{1,3,5\}, \{1,3,5\}, \{2,4,6\})$ is on Horn's 
list, as we have seen, the inequality
$$
	\gamma_2 + \gamma_4 + \gamma_6 \, \leq \, \alpha_1 + \alpha_3 + 
\alpha_5 + \beta_1 + \beta_3 + \beta_5 \tag12
$$
is redundant, and follows in fact from ($*$) and the fact that the eigenvalues 
are listed in decreasing order.  Indeed, let $\alpha_{\operatorname{ev}} = 
\alpha_2 + \alpha_4 + \alpha_6$, and $\alpha_{\operatorname{od}} = 
\alpha_1 + \alpha_3 + \alpha_5$, and similarly for $\beta$  and $\gamma$.  
The fact that each of $\alpha$, $\beta$, and $\gamma$  is nonincreasing 
implies the inequalities $\alpha_{\operatorname{ev}} \leq 
\alpha_{\operatorname{od}}$, $\beta_{\operatorname{ev}} \leq 
\beta_{\operatorname{od}}$, and $\gamma_{\operatorname{ev}} \leq 
\gamma_{\operatorname{od}}$.  But ($*$) says that 
$\gamma_{\operatorname{od}} + \gamma_{\operatorname{ev}} = 
\alpha_{\operatorname{od}} + \alpha_{\operatorname{ev}} + 
\beta_{\operatorname{od}} + \beta_{\operatorname{ev}}$.  Hence 
$\gamma_{\operatorname{ev}} \leq \alpha_{\operatorname{od}} + 
\beta_{\operatorname{od}}$, as asserted. 

In fact, for $n = 6$, this is the {\em only} triple on Horn's list that can be 
omitted (and the only triple that gives rise to an example like Example 1).  
As $n$ increases, however, the number of redundant triples on the list 
increases rapidly.  Describing the actual minimal set, however, requires 
notions from other areas of mathematics, and is postponed to later sections.

As this example indicates, one reason why the inequalities defining the 
realizable triples $(\alpha ,\beta ,\gamma)$ pose problems is because, 
together with equation ($*$) and the inequalities ($*_{IJK}$), one also has 
the $3n-3$ inequalities $\alpha_1\geq  \alpha_2\geq  \ldots  \geq  
\alpha_n$, $\beta_1\geq  \beta_2\geq  \ldots  \geq  \beta_n$, 
$\gamma_1\geq  \gamma_2\geq  \ldots  \geq  \gamma_n$.    

There is also a literature describing the possibilities for a particular 
eigenvalue $\gamma_k$ (with $k$ fixed) of $C = A + B$ in terms of the 
eigenvalues of $A$ and B.  In fact, this problem was solved much earlier (see 
[Joh]).  The result is that the $k^{\text{th}}$ largest eigenvalue of $C$ can 
take on any value in an interval:  
$$
	\operatornamewithlimits{Max}_{i+j = n+k} \alpha_i + \beta_j \, \leq \,
\gamma_k \, \leq \, \operatornamewithlimits{Min}_{i+j = k+1} \alpha_i + 
\beta_j. \tag13
$$
We have seen in (2) and (11) that $\gamma_k$ satisfies these inequalities, 
and since the set in question is a projection of a convex set in 
$\mathbb{R}^n$ (or the image of the connected space $U(n)$ by a 
continuous map) it must be an interval.  It therefore 
 suffices to produce Hermitian 
matrices to achieve each of the endpoints of the displayed interval.  Explicit 
diagonal matrices can be produced to achieve this [Joh].  This argument also 
shows that this question has the same answer for real symmetric matrices.  

More generally, one can specify a subset $K$ of $\{1, \ldots ,n\}$, and ask 
for the possible values of $\{\gamma_k \,|\, k \in K\}$, again with $\alpha$  
and $\beta$  given.  In principle, such results may be deduced from 
Theorem 1, although carrying this out does not appear to be easy.  One 
always has the inequalities ($*_{FGH}$) given for triples $(F,G,H)$, for $H$ 
a subset of $K$, and the duals (7) of inequalities ($*_{FGH}$) when $H$ 
contains the complement of $K$, together with 
 the inequalities saying that the 
eigenvalues in $K$ form a weakly decreasing sequence.  But already for $n = 
3$ and $K = \{1,3\}$ it is easy to see that such inequalities do not suffice: 
they do not imply that the missing $\gamma_2$, which is determined, is at 
least as large as $\gamma_3$. 

Not all inequalities that have been found involving eigenvalues of sums of 
Hermitian matrices are linear.  For example, M. Fiedler [Fi] showed that
$$
   \operatornamewithlimits{Min}_{\sigma \in S_n} \prod_{i=1}^n 
(\alpha_i+\beta_{\sigma(i)}) \, \leq \, \prod_{i=1}^n \gamma_i \, \leq\,  
\operatornamewithlimits{Max}_{\sigma \in S_n}
 \prod_{i=1}^n 
(\alpha_i+\beta_{\sigma(i)}). \tag14
$$

There are still some basic questions remaining.  One of these is the 
following:
\block
	Which triples $(I,J,K)$ of subsets of cardinality $r$ in $\{1, \ldots 
,n\}$ give true inequalities ($*_{IJK}$) for eigenvalues of all $n$ by $n$ 
Hermitian matrices $A$, $B$, C = A+B?  
\endblock
Let us call this set $H_r^n$.  This question was addressed by Horn [H2] and 
then by Zwahlen [Zw].  Theorem 1 asserts that $T_r^n \subset 
H_r^n$, but how much larger is $H_r^n$? 

For $r \leq 2$, the answer is quite simple, and was given by Horn and 
Zwahlen.  The set $H_1^n$ consists of triples with $i_1 + j_1 \leq k_1 + 
1$.  The set $H_2^n$ consists of triples such that $i_1 + j_1 \leq k_1 + 1$, 
$i_2 + j_1 \leq k_2 + 1$, $i_1 + j_2 \leq k_2 + 1$, and $i_1 + i_2 + j_1 + 
j_2 \leq k_1 + k_2 + 3$.  In fact, they show that for $r \leq 2$ and any 
triple not satisfying these conditions, there are diagonal matrices $A$, $B$, 
and $C = A+B$ whose eigenvalues violate the inequality ($*_{IJK}$).  They 
prove some partial results for $r = 3$.  On the basis of the evidence in these 
papers, it was natural to hope that, if $(I,J,K)$ is given to be in $U_r^n$, 
then it is in $H_r^n$ only if it is in $T_r^n$,  i.e., $H_r^n\cap U_r^n = 
T_r^n$.  However, this is not true:

\example{Example 2}  Take $r = 5$, $n = 25$, and $I = J = 
\{1,3,4,16,21\}$ and $K = \{5,10,15,20,25\}$.  This is in $H_r^n\cap 
U_r^n$.  Indeed, it is easy to see, by the same argument as in the proof of 
(12), that for any weakly decreasing $n$-tuples $\alpha$, $\beta$, and 
$\gamma$  such that $\sum \gamma_i = \sum \alpha_i + \sum \beta_i$, 
$$\align
   	\sum_{k \in K} \gamma_k = \sum_{s=0}^4 \gamma_{5+5s} &\leq 
	\sum_{s=0}^4 \alpha_{1+5s} + \sum_{s=0}^4 \beta_{1+5s} \\
	&\leq \sum_{i \in I} \alpha_i + \sum_{j \in J}  \beta_j.
\endalign$$
To see that $(I,J,K)$ is not in $T_5^{25}$, one can use the triple $(F,G,H)$ in 
$T_4^5$, with $F = G = \{1,2,4,5\}$ and $H = \{2,3,4,5\}$: $\sum_{f \in 
F} i_f +\sum_{g \in G} j_g = 82 > 80 =\sum_{h \in H} k_h + p(p+1)/2$.  
Many other examples can be constructed by this method, but none of them 
are very small. \endexample

Zwahlen [Zw] gave an example of a triple $(I,J,K)$ for which he produced 
violating matrices, but for which no diagonal matrices violate ($*_{IJK}$); 
this was, for $r = 3$, $n = 18$,  $I = J = \{1,6,11\}$, $K = \{2,9,18\}$.  
Thompson and Therianos [TT1] gave the simpler triple 
$(\{1,3,5\},\{1,3,5\},\{2,3,6\})$ for $r = 3$, $n = 6$, with the same 
property.  These triples, however, are not in $U_r^n$.  
 
\example{Example 3 (Buch)}  The triple $(I,J,K) = (\{1,3,5,6\}, \{1,3,5,6\}, 
\{2,3,6,9\})$ is in $U_4^9\smallsetminus T_4^9$, but no diagonal $9$ by 
$9$ Hermitian matrices can have eigenvalues violating ($*_{IJK}$).  If $A$ 
is the diagonal matrix with diagonal entries $(2,0,1,2,0,1,0,0,0)$, and $B$ 
is the direct sum of the matrix $\left(\smallmatrix 1/2 & \sqrt{3}/2 \\ 
\sqrt{3}/2 & 3/2 \endsmallmatrix \right)$ and the diagonal matrix with 
diagonal entries $(2,1,1,0,0,0,0)$, then the eigenvalues of $A$, $B$, and $C 
= A+B$ are $\alpha  = \beta  = (2,2,1,1,0,0,0,0,0)$, and $\gamma  = 
(3,3,3,1,1,1,0,0,0)$. These violate ($*_{IJK}$). \endexample

These examples indicate some of the subtlety of identifying triples that give 
correct inequalities, and finding violating matrices for triples that do not.

It is interesting to note that much of the extreme behavior can be detected 
by matrices that are 
diagonal or close to diagonal, at least in low dimensions.  
Indeed, this must have been how many of the inequalities were discovered.  
However, we will see that the proofs of the theorems are almost opposite to 
this: the eigenvectors of the matrices produced by the proofs are in general 
position.

The theorems of this section give a clear picture of the set of all 
triples $(\alpha,\beta,\gamma)$ that occur as spectra of Hermitian 
$n$ by $n$ matrices $A$, $B$, and $C = A + B$.  For small $n$ they give 
a reasonable set of inequalities that one can use to test if a particular 
triple $(\alpha,\beta,\gamma)$ occurs.  For large $n$, however, the number 
of inequalities to be tested increases dramatically (even if decreased 
by using Theorem 13 below).  When all of $\alpha$, $\beta$, and $\gamma$ 
are integral, one can test directly if $(\alpha,\beta,\gamma)$ occurs
by using Theorem 11 below. When the eigenvalues are rational, one can 
multiply them all by a common denominator to reduce to the integral case. 
We do not know a similarly direct criterion for arbitrary real eigenvalues.

The necessary conditions of the theorems extend readily from the realm of 
finite dimensional Hermitian operators to that of compact self-adjoint 
operators on a Hilbert space.  Indeed, it was in this context that Weyl [W] 
stated his results, cf. [Zw].  Such an operator $A$ has a sequence of 
positive eigenvalues $\alpha_1 \geq  \alpha_2 \geq  \alpha_3 \geq   
\ldots$    (each occurring according to multiplicity), and a similar sequence 
of negative eigenvalues.  For simplicity we consider only the positive 
eigenvalues: 

\proclaim{Theorem 6}  Suppose $A$, $B$, and $C = A+B$ are compact 
self-adjoint
operators on a Hilbert space, with $\alpha$, $\beta$, and $\gamma$ 
their sequences of positive eigenvalues, and assume that each of these 
sequences is infinite.  Then \rom{(}$*_{IJK}$\rom{)} holds for all $(I,J,K)$ 
in $T_r^n$, for all $n$ and $r$ with $r < n$. \endproclaim

The theorems here describe all possible eigenvalues of $A+B$ when $A$ 
and $B$ are $n$ by $n$ Hermitian matrices with eigenvalues $\alpha$ 
and $\beta$.  In contrast, P. Biane [Bi] has shown that, as $n$ gets 
large, for almost all choices of $A$ and $B$, the eigenvalues of $A+B$ 
are close to some $\gamma$ that depends only on $\alpha$ and $\beta$.

Much of the history of the eigenvalue problem before the recent events, 
together with 
an extensive bibliography, can be found in the survey [DST].

\head 2.  Invariant factors \endhead

We turn now to quite a different problem.  Consider an $n$ by $n$ matrix 
$A$ with coefficients in a discrete valuation ring $R$, whose determinant is 
not zero.  Let $\pi$  be a uniformizing parameter in $R$.  The cases 
that have been most studied are when $R = \mathbb{C}\{z\}$ is the ring of 
convergent power series in one variable, with $\pi   = z$, or when $R = 
\mathbb{Z}_p$ is the ring of $p$-adic integers, with $\pi   = p$.  By 
elementary row and column operations \footnote{The elementary row 
operations are: interchanging two rows; adding to any row a multiple by an 
element of $R$ times another row; multiplying any row by an invertible 
element of $R$. Similarly for elementary column operations.}, one can 
reduce $A$ to a diagonal matrix, with diagonal entries $\pi^{\alpha_1}, 
\pi^{\alpha_2}, \ldots , \pi^{\alpha_n}$, for unique nonnegative integers  
$\alpha_1 \geq  \ldots  \geq  \alpha_n$.  We call $\alpha  = (\alpha_1, 
\ldots , \alpha_n)$ the {\bf invariant factors} (or Smith invariants) of $A$.  
The question in this case is: 
\block
	Which $\alpha$, $\beta$, $\gamma$  can be the invariant factors of 
	matrices $A$, $B$, and $C$ if $C = A\cdot B$?
\endblock

In the case of convergent power series, this problem was proposed by I. 
Gohberg and M. A. Kaashoek, and was attacked particularly by Thompson 
and his coworkers (cf. [Th1], [Th2], [Th3], [Thi]). 

These matrices correspond to endomorphisms of $R^n$, with cokernels 
being torsion modules with at most $n$ generators.  Such a module is 
isomorphic to a direct sum 
$$
	R/\pi^{\alpha_1}\!R\, \oplus \, R/\pi^{\alpha_2} \! R \, \oplus 
\ldots \oplus \,
	R/\pi^{\alpha_n}\!R, \qquad  \alpha_1 \geq  \ldots  \geq  \alpha_n.
$$
We call $\alpha  = (\alpha_1, \ldots , \alpha_n)$ the {\bf invariant factors} 
of the module.  Denoting cokernels of $A$, $B$, and $C$ by $\mathcal{A}$, 
$\mathcal{B}$, and $\mathcal{C}$,  one has a short exact sequence
$$
	0 \to \mathcal{B} \to \mathcal{C} \to \mathcal{A} \to 0,
$$
i.e., $\mathcal{B}$ is a submodule of $\mathcal{C}$, with 
$\mathcal{C}/\mathcal{B} \cong \mathcal{A}$.  Conversely (cf. [Th3]), such 
an exact sequence corresponds to $n$ by $n$ matrices $A$, $B$, and $C$ 
with $ C = A\cdot B$.  This correspondence is seen easily by applying the 
``snake lemma'' (see e.g. [La] \S III.7) to the diagram
$$
\xymatrix{
	 R^n \ar[r]^B \ar[d]_C &  R^n \ar[d]^A \\
 	 R^n \ar[r]^{\operatorname{id}} & R^n     }
$$
Note that the equality ($*$), that $\sum \gamma_i = \sum \alpha_i + \sum 
\beta_i$, is obviously satisfied, since the determinant of a product is the 
product of the determinants.  
It is also not difficult to verify that inequality 
(1) must also be valid.  Thompson [Th1] proved (6) in this setting.

When $R$ is the ring of $p$-adic integers, one is asking what finite abelian 
$p$-groups can appear in a short exact sequence.  The aim of workers in 
this field was to use the numbers of submodules $\mathcal{B}$ of a given 
module $\mathcal{C}$ with specified invariant factors for $\mathcal{B}$, 
$\mathcal{C}$, and $\mathcal{A} = \mathcal{C}/\mathcal{B}$, as the 
structure constants to define a ring.  
In fact, this idea goes back to a lecture 
of E. Steinitz in 1900, but this was lost until the 1980's.  From the 1940's to 
the 1960's, this theory was developed by P. Hall, J.A. Green, and T. Klein. 
The conclusion of 
this study that is relevant to our question was that a triple 
$(\alpha,\beta,\gamma)$ occurs for $p$-groups, or in fact for any discrete 
valuation ring, if and only if a certain nonnegative integer 
$c_{\alpha\,\beta}^{\,\,\gamma}$, called the Littlewood-Richardson 
coefficient, is nonzero.  These coefficients, which have a purely 
combinatorial description, arose in representation theory, and will be 
defined and discussed in the next section.  
 
Although there seems to be little relation between this problem and the 
eigenvalue problem, they in fact have exactly the same answers.  In 
particular, this implies the fact that the answer to this problem is 
independent of the discrete valuation ring.

\proclaim{Theorem 7}   For any discrete valuation ring $R$, a triple 
$(\alpha,\beta,\gamma)$ occurs as the invariant factors of $A$, $B$, and $C 
= A\cdot B$ if and only if 
$\sum \gamma_i = \sum \alpha_i + \sum \beta_i$  and the inequalities 
\rom{(}$*_{IJK}$\rom{)} are satisfied for all $(I,J,K)$ in $T_r^n$, for all $r 
< n$. \endproclaim

In fact, these results can be extended to any principal ideal domain $R$, 
such as the integers $\mathbb{Z}$ or a polynomial ring $F[T]$ in one 
variable over a field.  In this case the invariant factors of a matrix $A$ with 
nonzero determinant are a chain $a_1 \subset a_2 \subset \ldots \subset 
a_n$ of nonzero ideals in $R$.  The matrix $A$ can be reduced by 
left and right multiplication by matrices 
in $GL(n,R)$ to a diagonal matrix whose successive 
entries generate the ideals in the chain; equivalently, the cokernel of $A$ is 
isomorphic to $\oplus_{i = 1}^n R/a_i$.  

\proclaim{Theorem 8}   For any principal ideal domain $R$, a triple 
$(a,b,c)$ of chains of ideals occurs as the invariant factors of $A$, $B$, and 
$C = A\cdot B$ (or of torsion modules $\mathcal{A}$, $\mathcal{B}$, and 
$\mathcal{C}$ with at most $n$ generators with $\mathcal{B} \subset 
\mathcal{C}$ and $\mathcal{C}/\mathcal{B} \cong \mathcal{A}$) if and only 
if $\prod_i c_i = \prod_i a_i \cdot \prod_i b_i$, and, for all $r < n$, and 
all $(I,J,K)$ in $T_r^n$, we have
$$
	\prod_{k \in K} c_k  \supset \prod_{i \in I} a_i \cdot  \prod_{j \in 
J} 	b_j.      
$$ \endproclaim

One can also consider matrices with entries in the quotient field of $R$, in 
which case the sequences $\alpha, \beta$, and $\gamma$  may contain 
negative integers, and the ideals $a$, $b$, and $c$ may include fractional 
ideals.  Theorems 7 and 8 extend immediately to these situations, by 
multiplying the matrices by scalars to get all entries in $R$.

This theorem also gives a solution of a problem called the {\em Carlson 
problem} ([C], cf. [Th2], [SQS]).  The general theorem is the following:

\proclaim{Theorem 9}  Let $A$ and $B$ be $p$ by $p$ and $q$ by $q$ 
matrices with entries in a principal ideal domain $R$, with nonzero 
determinants, and with invariant factors $a_1 \subset \ldots \subset a_p$ and 
$b_1 \subset \ldots \subset b_q$.  Then the possible invariant factors 
for a matrix of the form
$$
	C = \pmatrix A  & X \\ 0 & B \endpmatrix 
$$
for $X$ a $p$ by $q$ matrix with entries in $R$, 
are those $c_1
\subset \ldots \subset c_n$, with $n = p + q$, 
for which $\prod_{k=1}^n c_k = \prod_{i=1}^p a_i \cdot \prod_{j=1}^q 
b_j$, and 
$$
	\prod_{k \in K} c_k  \supset \prod \Sb
	i \in I \\ i \leq p \endSb  a_i \cdot  
\prod \Sb  j \in J \\ j \leq q \endSb  b_j
$$
for all $(I,J,K)$ in $T_r^n$ and all $r < n$. \endproclaim

An arbitrary $s$ by $s$ matrix $A$ with entries in a field $F$ has 
invariant factors $a_1 \subset \ldots \subset a_p$ in $F[T]$: these 
are the invariant factors of the $F[T]$-module $F^n$, where $T$ acts 
by $T\cdot v = A \, v$ for $v \in F^n$.
The following theorem solves the original Carlson problem ([C]):

\proclaim{Theorem 10}  Let $A$ and $B$ be $s$ by $s$ and $t$ by $t$ 
matrices over a field $F$.  Let $a_1 \subset \ldots \subset a_p$ and $b_1 
\subset \ldots \subset b_q$ be the invariant factors of $A$ and $B$.  Then 
the possible invariant factors of a matrix $C = \left(\smallmatrix A  & X \\ 0 
& B \endsmallmatrix\right)$, with $X$ an $s$ by $t$ matrix with entries in 
$F$, are those $c_1 \subset \ldots \subset c_n$, with $n = p+q$, that 
satisfy the conditions of Theorem 9. \endproclaim

For example, if $A = B = \left(\smallmatrix 0  & 1 \\ 0 & 0 
\endsmallmatrix\right)$, the invariant factors are $a_1 = b_1 = (T^2)$, 
with $p = q = 1$.  There are three possible invariant factors of $C$: $(T^2) 
\subset (T^2)$, which occurs with $X = \left(\smallmatrix 0  & 0 \\ 0 & 0 
\endsmallmatrix\right)$;  $(T^3) \subset (T)$, with $X = \left(\smallmatrix 
1  & 0 \\ 0 & 0 \endsmallmatrix\right)$;  and $(T^4) \subset (1)$, with $X 
= \left(\smallmatrix 0  & 0 \\ 1 & 0 \endsmallmatrix\right)$.  Note that 
$p$ and $q$ are usually smaller than $s$ and $t$, and the matrices $A$, 
$B$, and $X$ have different meanings in Theorems 9 and 10.  

For an application, suppose $A$ is a nilpotent matrix, with Jordan blocks of 
sizes $\alpha_1\geq  \alpha_2\geq  \ldots  \geq  \alpha_p$, and $B = 0$ is 
a $q$ by $q$ matrix.  We apply Theorem 10 with $a_i = t^{\alpha_i}$, $b_j = 
t^{\beta_j}$ with $\beta_j = 1$, and $c_k = t^{\gamma_k}$.  In this case a 
rule of Pieri, given in equation (16) in the next section, implies that the 
possible $C$ have Jordan blocks of sizes $\gamma_1 \geq  \ldots  \geq  
\gamma_t$, with $p \leq t \leq p+q$, $\sum \gamma_i= \sum \alpha_i + 
p$, $\alpha_i+1 \geq  \gamma_i \geq  \alpha_i$ for $1 \leq i \leq p$, and 
$\gamma_i = 1$ for $p < i \leq t$.   This result has been used by C. R. 
Johnson and E. A. Schreiner [JS] to give a quick proof of a theorem of H. 
Flanders characterizing which pairs $(C,D)$ of an $m$ by $m$ matrix $C$ 
and an $n$ by $n$ matrix $D$ have the form $(A\cdot B,B\cdot A)$ for 
some $m$ by $n$ matrix $A$ and some $n$ by $m$ matrix B.

Unlike the situations we have seen previously, in Theorems 9 and 10 the 
matrices $A$ and $B$, and not merely their invariant factors, can be 
specified arbitrarily in advance.	

The inequalities (13) also tell, for a fixed $k$ between $1$ and $n$, and 
fixed partitions $\alpha$  and $\beta$, exactly which integers can be 
$\gamma_k$ for some triple $(\alpha,\beta ,\gamma)$ that occur as 
invariant factors.  As before, the necessity we know. To prove sufficiency 
it suffices to construct matrices that realize the possibilities.  This has 
been done by J. F. Queir\'o and E. Marques de S\'a [QS].  
In light of the stronger 
results we now have, this raises the question of what can be said about other 
subsets.

\head 3.  Highest weights \endhead

An irreducible, finite-dimensional, holomorphic representation of 
${GL}_n(\mathbb{C})$ is characterized by its highest weight, which is a 
weakly decreasing sequence 
$$
	\alpha  = (\alpha_1 \geq  \alpha_2 \geq  \ldots  \geq  \alpha_n)
$$
of integers.  For example, the representation $\bigwedge^k(\mathbb{C}^n)$ 
corresponds to the sequence $(1,1, \ldots , 1, 0, \ldots , 0)$ consisting of 
$k$ $1$'s and $n-k$ $0$'s, and the representation 
$\operatorname{Sym}^k(\mathbb{C}^n)$ has highest weight $(k, 0, \ldots , 
0)$.  In general such a representation contains a nonzero vector $v$ (called 
a highest weight vector), such that for any upper triangular matrix $X$ in 
${GL}_n(\mathbb{C})$, whose entries down the diagonal are $x_1, x_2, 
\ldots , x_n$,
 $$
X \cdot v = x_1^{\alpha_1} \, x_2^{\alpha_2} \cdots x_n^{\alpha_n} 
\,\, v.
$$
We denote the irreducible representation with highest weight $\alpha$  by  
$V(\alpha)$.
It is a basic fact of representation theory that 
${GL}_n(\mathbb{C})$ is {\em reductive}. This means that any 
finite-dimensional holomorphic representation decomposes into a direct 
sum of irreducible representations, and the number of times that a given 
irreducible representation $V(\gamma)$ appears 
in the sum is independent of the choice 
of the decomposition.  In particular, for any $\alpha$ and $\beta$, the 
tensor product $V(\alpha) \otimes V(\beta)$ decomposes into a direct sum 
of representations $V(\gamma)$.  Define $c_{\alpha\,\beta}^{\,\,\gamma}$  
to be the number of copies of $V(\gamma)$ in an irreducible decomposition 
of $V(\alpha) \otimes V(\beta)$.  The problem of interest in this situation is: 
\block
	When does $V(\gamma)$ occur in $V(\alpha) \otimes V(\beta)$, i.e., 
	when is $c_{\alpha\,\beta}^{\,\,\gamma}  > 0$?
\endblock
It follows immediately from the definition of highest weights that a 
necessary condition for this is that $\sum \gamma_i = \sum \alpha_i + 
\sum \beta_i$.  Other conditions are more difficult to find, although an 
expert may attempt to prove some of the inequalities (1) -- (7).

A simple case of this problem 
is when $\beta  = (1, \ldots , 1)$, so $V(\beta)$ is the 
one-dimensional determinant representation.  In this case $V(\alpha) 
\otimes V(\beta)$ is equal to 
$V(\alpha_1+1, \alpha_2+1, \ldots , \alpha_n+1)$.  In 
particular, the problem is unchanged if each of the representations is 
tensored by this determinant representation several times.  Therefore we 
may assume that each of $\alpha$, $\beta$, and $\gamma$  consists of 
nonnegative integers, i.e., is a partition.  Equivalently, one need 
only consider 
{\em polynomial} representations of ${GL}_n(\mathbb{C})$.  These 
representations were constructed in the nineteenth century by J. Deruyts [D].
For a simple construction of them, see [Fu2], \S8.1.

In 1934 Littlewood and Richardson [LR] gave a remarkable combinatorial 
formula for the numbers $c_{\alpha\,\beta}^{\,\,\gamma}$, although a 
complete proof of the formula 
was only given much later (cf. [Mac], [Fu2]).  This formula is 
stated in terms of the Young (or Ferrers) diagrams corresponding to the 
partitions.  The diagram of $\alpha$ is an array of boxes, lined up 
at the left, 
with $\alpha_i$ boxes in the $i^{\text{th}}$ row, with the rows 
arranged from top to 
bottom.  For example, 
$$
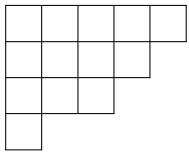
\grub0.4pt
\hbox{\kern-30\jedn
\MT 0 0 \LT 0 24
\MT 6 0 \LT 6 24
\MT 12 6 \LT 12 24
\MT 18 6 \LT 18 24
\MT 24 12 \LT 24 24
\MT 30 18 \LT 30 24
\MT 0 0 \LT 6 0
\MT 0 6 \LT 18 6
\MT 0 12 \LT 24 12
\MT 0 18 \LT 30 18
\MT 0 24 \LT 30 24
}
$$
is the Young diagram of $(5,4,3,1)$.  We follow the convention, reinforced 
by the Young diagrams, of identifying two partitions if they differ by a string 
of zeros at the end.

Their rule states that $c_{\alpha\,\beta}^{\,\,\gamma}$  is zero unless 
$\sum \gamma_i = \sum \alpha_i + \sum \beta_i$ and 
the Young diagram of $\alpha$  is contained in that for $\gamma$, i.e., 
$\alpha_i \leq \gamma_i$ for $1 \leq i \leq n$.  The complement of the 
Young diagram of $\alpha$ in that of $\gamma$, denoted $\gamma 
\smallsetminus \alpha$, then consists of $\sum \beta_i$
boxes.  We {\bf order} the boxes of $\gamma 
\smallsetminus \alpha$  by first listing the boxes in the top row, from right 
to left, then the boxes in the second row from right to left, and so on down 
the array. The {\bf Littlewood-Richardson coefficient} 
$c_{\alpha\,\beta}^{\,\,\gamma}$ is the number of ways to fill the boxes of 
$\gamma \smallsetminus \alpha$, with one integer in each box, so that the 
following four conditions are satisfied:
\roster
\item"(i)"  	The entries in any row are weakly increasing from left to right;
\item"(ii)"  	The entries in each column are strictly increasing from top to 
bottom.
\item"(iii)"  The integer $i$ occurs exactly $\beta_i$ times. 
\item"(iv)"  For any $p$ with $1 \leq p < \sum \beta_i$, and any $i$ with 
$1 \leq i < n$, the number of times $i$ occurs in the first $p$ boxes of the 
ordering is at least as large as the number of times that $i+1$ occurs in 
these first $p$ boxes.
\endroster
(The last condition says that when the 
entries are listed, from right to left in 
rows, from the top row to the bottom, they form a {\em lattice word.})  

For $\alpha  = (3,2,1)$, $\beta  = (3,2,2)$, and $\gamma  = (5,4,3,1)$, the 
following are some of the ways to fill the boxes of $\gamma \smallsetminus 
\alpha$  satisfying the first three conditions.  
$$
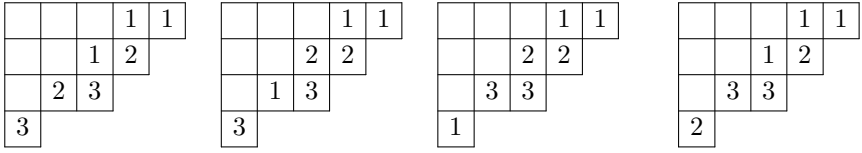
\grub0.4pt
\hbox{\kern-30\jedn
\MT 0 0 \LT 0 24
\MT 6 0 \LT 6 24
\MT 12 6 \LT 12 24
\MT 18 6 \LT 18 24
\MT 24 12 \LT 24 24
\MT 30 18 \LT 30 24
\MT 0 0 \LT 6 0
\MT 0 6 \LT 18 6
\MT 0 12 \LT 24 12
\MT 0 18 \LT 30 18
\MT 0 24 \LT 30 24
\cput(3,2,3) \cput(9,8,2) \cput(15,8,3)
\cput(15,14,1) \cput(21,14,2) \cput(21,20,1) \cput(27,20,1)
\kern36\jedn
\MT 0 0 \LT 0 24
\MT 6 0 \LT 6 24
\MT 12 6 \LT 12 24
\MT 18 6 \LT 18 24
\MT 24 12 \LT 24 24
\MT 30 18 \LT 30 24
\MT 0 0 \LT 6 0
\MT 0 6 \LT 18 6
\MT 0 12 \LT 24 12
\MT 0 18 \LT 30 18
\MT 0 24 \LT 30 24
\cput(3,2,3) \cput(9,8,1) \cput(15,8,3)
\cput(15,14,2) \cput(21,14,2) \cput(21,20,1) \cput(27,20,1)
\kern36\jedn
\MT 0 0 \LT 0 24
\MT 6 0 \LT 6 24
\MT 12 6 \LT 12 24
\MT 18 6 \LT 18 24
\MT 24 12 \LT 24 24
\MT 30 18 \LT 30 24
\MT 0 0 \LT 6 0
\MT 0 6 \LT 18 6
\MT 0 12 \LT 24 12
\MT 0 18 \LT 30 18
\MT 0 24 \LT 30 24
\cput(3,2,1) \cput(9,8,3) \cput(15,8,3)
\cput(15,14,2) \cput(21,14,2) \cput(21,20,1) \cput(27,20,1)
\kern40\jedn
\MT 0 0 \LT 0 24
\MT 6 0 \LT 6 24
\MT 12 6 \LT 12 24
\MT 18 6 \LT 18 24
\MT 24 12 \LT 24 24
\MT 30 18 \LT 30 24
\MT 0 0 \LT 6 0
\MT 0 6 \LT 18 6
\MT 0 12 \LT 24 12
\MT 0 18 \LT 30 18
\MT 0 24 \LT 30 24
\cput(3,2,2) \cput(9,8,3) \cput(15,8,3)
\cput(15,14,1) \cput(21,14,2) \cput(21,20,1) \cput(27,20,1)
}
$$
The first three examples satisfy the fourth condition. The fourth example 
does not, since the first six boxes in the ordering have more $3$'s than 
$2$'s.  One sees easily that the first three are the only possibilities satisfying 
all four conditions, so $c_{\alpha\,\beta}^{\,\,\gamma}  = 3$.  

Two special cases of this rule were known to Pieri in the context of 
Schubert calculus.  Let $\alpha  = (\alpha_1, \ldots , \alpha_n)$.  If $\beta  
= (p)$, then the possible $\gamma$  for which 
$c_{\alpha\,\beta}^{\,\,\gamma}  \ne 0$ are those of the form 
$(\gamma_1, \ldots , \gamma_{n+1})$ with 
$$\aligned
\gamma_1 \geq  \alpha_1 \geq  \gamma_2 \geq  \alpha_2 \geq  \ldots  
\geq  \gamma_n \geq  &\alpha_n \geq  \gamma_{n+1} \geq  0, \\ 
&\text{ with } \sum \gamma_i = \sum \alpha_i + p.
\endaligned \tag15$$
In these cases $c_{\alpha\,\beta}^{\,\,\gamma}  = 1$. (For 
representations of $GL_n(\mathbb{C})$ only those with $\gamma_{n+1} = 
0$ are allowed.)  In terms of Young diagrams, $\beta$  consists of a row of 
$p$ boxes, and the diagram of $\gamma$  is obtained from that of $\alpha$  
by adding $p$ boxes, with no two in any column.  

The other Pieri rule is for 
$\beta  = (1, \ldots ,1)$ consisting of $p$ $1$'s. The possible 
$\gamma$  with $c_{\alpha\,\beta}^{\,\,\gamma}  \ne 0$ also have 
$c_{\alpha\,\beta}^{\,\,\gamma}  = 1$, and these have the form  $\gamma  
= (\gamma_1, \ldots , \gamma_t)$  with
$$
	\alpha_i + 1 \geq  \gamma_i \geq  \alpha_i \,\,\, \text{ for all $i$, 
and }\,\,\,  \sum \gamma_i = \sum \alpha_i + p. \tag16
$$
Here $\beta$  is a column of $p$ boxes, and the diagram of $\gamma$  is 
obtained from that of $\alpha$  by adding $p$ boxes, with no two in any 
row.

Another case that is easy to analyze directly from the Littlewood-Richardson 
rule is the case when $\gamma_i = \alpha_i + \beta_i$ for all $i$. In this 
case one sees that $c_{\alpha\,\beta}^{\,\,\gamma}  = 1$.  These are the 
partitions which correspond to the triples of subsets listed in (6).

\proclaim{Theorem 11}  The Littlewood-Richardson coefficient 
$c_{\alpha\,\beta}^{\,\,\gamma}$  is positive exactly when $\sum 
\gamma_i = \sum \alpha_i + \sum \beta_i$ and the inequalities 
\rom{(}$*_{IJK}$\rom{)} are valid for all $(I,J,K)$ in $T_r^n$, and all $r < 
n$. \endproclaim

Equivalently, in light of Theorem 1, $c_{\alpha\,\beta}^{\,\,\gamma}$ is 
positive exactly when there is a triple of $n$ by $n$ Hermitian matrices 
$A$, $B$, $C = A+B$ with eigenvalues $\alpha$, $\beta$, and $\gamma$.
For example, if there are permutations $\sigma$ and $\tau$ in $S_n$ such 
that $\alpha_{\sigma(i)} + \beta_{\tau(i)} = \gamma_i$ for $1 \leq i \leq 
n$, then $c_{\alpha\,\beta}^{\,\,\gamma}$ must be positive, since there 
are diagonal matrices $A$, $B$, and $C$ with these eigenvalues and 
$C=A+B$.  For a combinatorial proof of this fact, which is a special case of 
the general (proven) ``PRV conjecture,'' see [KT].

The first triple with $c_{\alpha\,\beta}^{\,\,\gamma}$  greater than $1$, 
i.e., with the smallest $\sum \gamma_i$,  has $\alpha  = (2,1,0)$, $\beta  = 
(2,1,0)$, and $\gamma  = (3,2,1)$.  In this case the coefficient is $2$, 
because there are two arrays that satisfy (i) -- (iv):
$$
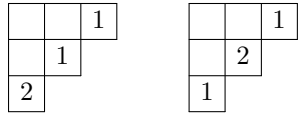
\grub0.4pt
\hbox{\kern-18\jedn
\MT 0 0 \LT 0 18
\MT 6 0 \LT 6 18
\MT 12 6 \LT 12 18
\MT 18 12 \LT 18 18
\MT 0 0 \LT 6 0
\MT 0 6 \LT 12 6
\MT 0 12 \LT 18 12
\MT 0 18 \LT 18 18
\cput(3,2,2)\cput(9,8,1)\cput(15,14,1)
\kern30\jedn
\MT 0 0 \LT 0 18
\MT 6 0 \LT 6 18
\MT 12 6 \LT 12 18
\MT 18 12 \LT 18 18
\MT 0 0 \LT 6 0
\MT 0 6 \LT 12 6
\MT 0 12 \LT 18 12
\MT 0 18 \LT 18 18
\cput(3,2,1)\cput(9,8,2)\cput(15,14,1) }
$$
Note that this triple corresponds to the triple $(\{1,3,5\}, \{1,3,5\}, 
\{2,4,6\})$ of subsets, which appeared in examples in Section 1.  The fact 
that this was the first triple to give a redundant Horn condition is not a 
coincidence.  In fact, let us define $R_r^n$ to be the set of triples $(I,J,K)$ 
of subsets of cardinality $r$ in $\{1, \ldots , n\}$ such that the 
corresponding triple $(\lambda,\mu,\nu) = 
(\lambda(I),\lambda(J),\lambda(K))$ has $c_{\lambda \,\mu}^{\,\,\nu} = 
1$:
$$
	R_r^n = \{(I,J,K) \in U_r^n \,\,| \,\, c_{\lambda \, \mu}^{\,\,\nu} = 
1\}. \tag17
$$
Similarly, let
$$
	S_r^n = \{(I,J,K) \in U_r^n \,\,|\,\, c_{\lambda \, \mu}^{\,\,\nu} 
\ne 0\}. \tag18
$$

The following theorem says that the triples occuring in Horn's conjecture are 
determined by the Littlewood-Richardson rule.

\proclaim{Theorem 12}  For all $r < n$, the sets $S_r^n$ and $T_r^n$ are 
equal. \endproclaim

The sets of triples therefore have the following inclusions:
$$
	  R_r^n \subset S_r^n  = T_r^n \subset U_r^n. \tag19
$$
\proclaim{Theorem 13}  Given that $\sum \gamma_i = \sum \alpha_i + 
\sum \beta_i$, and that each of $\alpha$, $\beta$, and $\gamma$ forms a 
weakly decreasing sequence, the conditions \rom{(}$*_{IJK}$\rom{)} for all 
$(I,J,K)$ in $T_r^n$ and all $r < n$, are implied by the conditions 
\rom{(}$*_{IJK}$\rom{)} for all $(I,J,K)$ in $R_r^n$ and all $r < n$. 
\endproclaim

In the theorems involving inequalities ($*_{IJK}$), Horn's sets $T_r^n$ can 
therefore be replaced by the smaller sets $R_r^n$. (An indication of how 
much smaller these sets are can be found at the end of Section 9.) 
 In fact, Knutson, Tao, 
and Woodward  
have announced [KTW] that the conditions ($*_{IJK}$) for $(I,J,K)$ in $R_r^n$  
are independent, i.e., none of them can be omitted.  

The Littlewood-Richardson coefficients also arise from the representation 
theory of finite symmetric groups.  The complex irreducible representations 
of the symmetric group $S_a$ are indexed by partitions of $a$; for a 
partition $\alpha$  of $a$, let $V_\alpha$  be this representation.  If  
$\alpha$  is a partition of $a$, and $\beta$  is a partition of $b$, then the 
representation $V_\alpha \otimes V_\beta$  is a representation of  $S_a 
\times  S_b$.  The Littlewood-Richardson number 
$c_{\alpha\,\beta}^{\,\,\gamma}$  is the number of times $V_\gamma$  
appears in the representation of $S_{a+b}$ that is induced from $V_\alpha 
\otimes V_\beta$  by the standard inclusion of $S_a \times S_b$ in 
$S_{a+b}$ (cf. [Mac], \S I.7 or [Fu2], \S7.3).  Theorems 11 and 12 therefore  
characterize the representations that occur in this decomposition.

If $\alpha$ and $\beta$ are given, the interval in (13) specifies exactly what 
can be the length $\gamma_k$ of the $k^{\text{th}}$ row of those 
$\gamma$ for which $c_{\alpha\,\beta}^{\,\,\gamma}$ is not zero; this 
follows from the equivalent problem discussed in Section 2.  Again there is 
the interesting open question of specifying the possible lengths 
$\{\gamma_k \,|\, k \in K\}$ of a prescribed subset of rows of such 
$\gamma$.

\head 4.  Schubert calculus \endhead

Let $X = Gr(n,\mathbb{C}^m)$ be the Grassmann variety of $n$-dimensional 
subspaces $L$ of $\mathbb{C}^m$.  It is a complex manifold of complex 
dimension $n(m-n)$.  A {\bf complete flag} $F_{\sssize{\bullet}}$ is a 
chain $ 0 = F_0 \subset F_1 
\subset F_2 \subset \ldots \subset F_m = \mathbb{C}^m$ of subspaces of 
$\mathbb{C}^m$, with $\operatorname{dim}(F_i) = i$ for all $i$.
For any subset $P = \{p_1 < p_2 < \ldots < p_n\}$ of 
cardinality $n$ in $\{1, \ldots , m\}$, there is a {\bf Schubert variety} 
$\Omega_P(F_{\sssize{\bullet}})$ in $X$ defined by 
$$
	   \Omega_P(F_{\sssize{\bullet}}) = \{ L \in X \, |\, 
\operatorname{dim}(L\cap F_{p_i}) \geq  i \,\, \text{ for } \,\,
1 \leq i \leq n \}. 
\tag20
$$
This is an irreducible closed subvariety of $X$ of dimension $\sum (p_i - 
i)$.  The homology classes $\omega_P = [\Omega_P(F_{\sssize{\bullet}})]$ 
of these varieties are 
independent of choice of the flag, and they form a basis for the integral 
homology of $X$.  For each partition $\alpha$  with  
$$
	m-n \geq  \alpha_1 \geq  \alpha_2 \geq  \ldots  \geq  \alpha_n \geq  
	0, 
$$
let $|\alpha| = \sum \alpha_i$, and define $\sigma_\alpha$ to be the 
cohomology class in $H^{2|\alpha|}X$ whose cap product with the 
fundamental class of $X$ is the homology class 
$\omega_P$, where $P$ is defined by setting $p_i = 
m - n + i - \alpha_i$.  Identifying cohomology and homology by 
this Poincar\'e duality isomorphism, we have:    
$$
		\sigma_\alpha  = \omega_P, \,\, \text{ where }\,\, \alpha_i = 
		m - n + i - p_i, \,\,\, 1 \leq i \leq n. \tag21
$$

These classes $\sigma_\alpha$  form a $\mathbb{Z}$-basis for the 
cohomology ring.  It follows that for any such partitions $\alpha$ and 
$\beta$, there is a unique expression
$$
	\sigma_\alpha \cdot \sigma_\beta  = \sum  d_{\alpha \, 
	\beta}^{\,\,\gamma} \,\,  \sigma_\gamma ,
$$
for integers $d_{\alpha \, \beta}^{\,\,\gamma}$, the sum over all 
$\gamma$  with $\sum \gamma_i = \sum \alpha_i + \sum \beta_i$.  It is a 
consequence of the fact that ${GL}_m(\mathbb{C})$ acts transitively on $X$ 
that all these coefficients are nonnegative.  The problem in this context is:
\block
	When does $\sigma_\gamma$  appear in the product $\sigma_\alpha 
	\cdot \sigma_\beta$, i.e., when is the coefficient $d_{\alpha \, 
	\beta}^{\,\,\gamma}$ positive?
\endblock

Algebraic geometers in the nineteenth century, especially Schubert, Pieri, 
and Giambelli, gave algorithms for computing in this cohomology ring (in 
spite of the fact that cohomology was not invented until many decades later).  
The Giambelli formula is a determinantal expression for any 
$\sigma_\alpha$  in terms of the basic classes $\sigma_k$, $1 \leq k \leq 
m-n$, where  $\sigma_k = \sigma_{(k)}$, and $(k)$ denotes the partition 
$(k, 0, \ldots , 0)$.  
The (original) Pieri formula (see (15)) gives the product 
$\sigma_k\cdot \sigma_\beta$.  However, the classical geometers did not 
give a general closed formula for the coefficients $d_{\alpha \, 
\beta}^{\,\,\gamma}$, nor did they give a criterion for these coefficients to 
be positive.

\proclaim{Theorem 14}  The class $\sigma_\gamma$  occurs in 
$\sigma_\alpha \cdot \sigma_\beta$  exactly when $\sum \gamma_i = 
\sum \alpha_i + \sum \beta_i$ and the inequalities 
\rom{(}$*_{IJK}$\rom{)} are valid for all $(I,J,K)$ in $T_r^n$, and all $r < 
n$. \endproclaim

This is stated for complex varieties and usual cohomology, but the same is 
true over any field, if one uses Chow groups instead of homology.  

For later use we state another of the fundamental facts about intersection 
theory on Grassmannians (cf. [Fu2], \S9.3).  Note that 
$$
	H^{2n(m-n)}X = \mathbb{Z}\cdot \sigma_\rho, \tag22
$$
where $\rho$ consists of the integer $m-n$ repeated $n$ times, i.e., 
$\sigma_\rho = \omega_{\{1,\ldots,r\}}$ is the class of a point in $X$. 
We identify this top 
cohomology group with $\mathbb{Z}$.  The dual of a class  
$\sigma_\gamma$  is the class $\sigma_{\gamma^\prime}$,  where  
$\gamma^\prime = (m-n-\gamma_n, \ldots , m-n-\gamma_2, m-n-
\gamma_1)$.  This means that, for $|\delta| + |\gamma| = n(m-n)$, the 
intersection number $\sigma_\gamma \cdot \sigma_\delta$ in $H^{2n(m-
n)}X = \mathbb{Z}$ is $1$ if $\delta = \gamma^\prime$, and $0$ 
otherwise.  Therefore, when $\sum \gamma_i = \sum \alpha_i + 
\sum \beta_i$, the class $\sigma_\gamma$  occurs in $\sigma_\alpha 
\cdot \sigma_\beta$  exactly when the product $\sigma_\alpha \cdot 
\sigma_\beta \cdot \sigma_{\gamma^\prime}$ is not zero.  In fact, 
$$
	\sigma_\alpha \cdot \sigma_\beta \cdot \sigma_{\gamma^\prime} = 
	d_{\alpha \, \beta}^{\,\,\gamma} \, \sigma_\rho. \tag23
$$

We will use an important fact that is part of this Schubert calculus.  Suppose, 
for each $1 \leq s \leq k$, $P(s) = \{p_1(s) < \ldots < p_n(s)\}$ is a subset 
of $\{1, \ldots , m\}$ of cardinality $n$, and $F_{\sssize{\bullet}}(s)$ is any 
complete flag.  Then the intersection of the Schubert varieties $\bigcap_ 
{s=1}^k \Omega_{P(s)}(F_{\sssize{\bullet}}(s))$ must be nonempty if the 
corresponding product $\prod_{s=1}^k \omega_{P(s)}$ is not zero.  This is 
a special case of a general fact in intersection theory, that the intersection of 
classes of varieties has a representative on the intersection of the varieties 
(cf. [Fu1]).  The converse is also true if 
the flags $F_{\sssize{\bullet}}(s)$ are 
in general position.  That is, there is a dense open set (the complement of a 
closed algebraic subset in the product of $k$ flag varieties) of $k$-tuples of 
flags which are in general position.  This is a special case of Kleiman's 
transversality theorem [Kle].

\head 5.  Singular values of sums and products \endhead

We regard $\mathbb{C}^m$ as the space of column vectors, and we use the 
standard Hermitian inner product, denoted $(u,v) = \sum_{i=1}^n u_i 
\overline{v_i}$ for $u$ and $v$ in $\mathbb{C}^m$; $|u|$ denotes the 
square root of $(u,u)$. 

An arbitrary real or complex $m$ by $n$ matrix $A$ has {\bf singular 
values} $a_1 \geq  a_2 \geq  \ldots  \geq  a_q \geq  0$, where $q$ is the 
minimum of $m$ and $n$.  They can be defined as follows.  The largest 
singular value $a_1$ is the maximum of $|(Au_1,v_1)|$ as $ u_1$ varies over 
unit vectors in $\mathbb{C}^n$, and $v_1$ varies over unit vectors in 
$\mathbb{C}^m$. Choosing $u_1$ and $v_1$ to achieve this maximum, define 
$$
	a_2 = \operatornamewithlimits{Max} \Sb |u_2| = 1, 
	\, u_2 \perp u_1 \\ |v_2| = 1, \, v_2 \perp v_1 \endSb \,\, |(A\, 
u_2,v_2)|.
$$
And so on for $a_3, \ldots , a_q$.  In fact, there is a unitary $n$ by $n$ 
matrix $U$, a unitary $m$ by $m$ matrix $V$, and an $m$ by $n$ 
matrix $D$ that has $a_1, \ldots , a_q$ down the diagonal and is otherwise 
zero, such that $A = V\, D\, U^*$.  This presentation 
is called a {\bf singular value 
decomposition} of $A$.  The vectors $u_i$ and $v_i$ can be taken as the 
$i^{\text{th}}$ columns of $U$ and $V$ respectively, for $1 \leq i \leq q$.  
This decomposition is important in multivariate analysis, where the entry 
$a_{ij}$ describe the $j^{\text{th}}$ property of an $i^{\text{th}}$ object; 
changing bases by $U$ and $V$ gives the best coordinates to distinguish 
correlations among the properties.  It is also used in numerical algorithms.  
Geometrically, if $A$ is real, the image $A(S^{n-1})$ of the unit sphere is an 
ellipsoid in $\mathbb{R}^m$, and the lengths of its principal axes are twice 
the singular values of $A$.

These singular values can also be determined in other ways.  For example, if 
$m \leq n$, then $a_1^2, \ldots , a_m^2$ are the eigenvalues of the 
Hermitian matrix $A\, A^*$; if $n \leq m$, one can similarly use $A^* A$.  
In general, $\pm A$ and $\pm A^*$ have the same singular values as $A$.  
A characterization we will use is that the eigenvalues of the $m+n$ by 
$m+n$ Hermitian matrix
 $$ 
\pmatrix 0 & A \\ A^* & 0 \endpmatrix 
$$
are 
$
a_1 \geq  a_1 \geq  \ldots  \geq  a_q \geq  0 \geq  \ldots  \geq  0 
\geq  -a_q \geq  \ldots  \geq  -a_2 \geq  -a_1
$,
where the number of zeros in the middle is the difference between $m$ and 
$n$.  The problem we consider in this setting is:
\block
What $a$, $b$, and $c$ can be the singular values of $m$ by $n$ 
matrices $A$, $B$, and $C$, with $C = A + B$?
\endblock

Given $a = (a_1 \geq  \ldots  \geq  a_q \geq  0)$, set 
$$
	\alpha  = (a_1 \geq  \ldots  \geq  a_q \geq  0 \geq  \ldots  \geq  
	0 \geq  -a_q \geq  \ldots  \geq  -a_1), 
$$
a nonincreasing sequence of length $m+n$, and define similarly sequences 
$\beta$  and $\gamma$  from $b$ and $c$.  For each triple $(I,J,K)$ of 
subsets of $\{1, \ldots , m+n\}$ of the same cardinality $r < m+n$, the 
inequality ($*_{IJK}$) on $\alpha$, $\beta$, and $\gamma$  determines a 
corresponding inequality on $a$, $b$, and $c$, that we denote by 
($**_{IJK}$).  For example, if $m = n = 2$, $r = 1$, and $I = \{3\}$, $J = 
\{1\}$, $K = \{3\}$, the inequality $\gamma_3 \leq \alpha_3 + \beta_1$ 
becomes the inequality $-c_2 \leq -a_2 + b_1$, or $a_2 \leq c_2 + b_1$.
In general, setting $I^\prime = \{i \in \{1, \dots , m+n\} \, | \, 
m+n+1-i \in I \}$, and similarly for $J$ and $K$, the inequality is 
$$
\sum \Sb k \in K \\ k \leq q \endSb c_k \, - \, 
\sum \Sb k \in K^\prime \\ k \leq q \endSb c_k \,\,\,\, \leq \,\,\,\, 
\sum \Sb i \in I \\ i \leq q \endSb a_i \, - \, 
\sum \Sb i \in I^\prime \\ i \leq q \endSb a_i \, + \, 
\sum \Sb j \in J \\ j \leq q \endSb b_j \, - \, 
\sum \Sb j \in J^\prime \\ j \leq q \endSb b_j. \tag$**_{IJK}$ 
$$
 
\proclaim{Theorem 15}  A triple $a$, $b$, $c$ occurs as singular values of 
$m$ by $n$ complex matrices $A$, $B$, and $C$ with $C = A+B$, if and 
only if the inequalities \rom{(}$**_{IJK}$\rom{)} are satisfied for all 
$(I,J,K)$ in $T_r^{m+n}$, for all $r < m+n$. \endproclaim

Many special cases of the necessity of these conditions had been known, 
cf. [AM], [TT1], [QS].

In contrast to the eigenvalue problem discussed in the first paragraph, the 
situation for real matrices is {\em not} the same as that for complex 
matrices.  For $1$ by $1$ complex matrices, a triple of nonnegative 
numbers occurs whenever each is at most the sum of the other two; for 
real matrices, one must be {\em equal} to the sum of the other two.  Here is 
a more interesting example:  

\example{Example 4}  For $3$ by $3$ matrices, the triple $a = (1,1,0)$, $b 
= (1,1,0)$, and $c = (1,1,1)$ occurs as singular values of complex matrices 
$A$, $B$ and $C = A + B$, but they do not occur for real matrices.  In the 
complex case, $A$ and $B$ can be diagonal matrices with entries $(1, \zeta, 
0)$ and $(0, \zeta^{-1}, 1)$, with $\zeta$ a primitive $6^{\text{th}}$ root 
of unity.  It is a straightforward calculation to verify that they do not occur 
for real matrices.  \endexample

When one writes out the inequalities of Theorem 15, one finds that most of 
them are redundant, and that the essential inequalities remaining have quite 
a simple form.  For example, for $4$ by $4$ matrices, the triple 
$(\{1,3,5\},\{2,3,5\},\{2,4,7\})$ in $T_3^8$ gives the 
inequality $c_2 + c_4 - c_2 \leq a_1 + a_3 - a_4 + b_2 + b_3 - b_4$. 
This inequality  
follows from the inequality $c_4 \leq a_3 + b_2$ and the inequalities 
$a_4 \leq a_1$ and $b_4 \leq b_3$. 
The triple $(\{3,7\}, \{2,3\}, \{4,8\})$ in 
$T_2^8$, on the other hand, gives the inequality $c_4 - c_1 \leq a_3 - a_2 
+ b_2 + b_3$, or 
$$
	a_2 + c_4 \leq c_1 + a_3 + b_2 + b_3, \tag24
$$
which is in fact essential: it determines one of the facets of the cone defined 
in Theorem 15.  

For simplicity, let us consider only $n$ by $n$ matrices.  For $n = 1$, 
minimal inequalities defining the possible singular values are:  $c_1 \leq a_1 
+ b_1$,  $b_1 \leq a_1 + c_1$, and $a_1 \leq b_1 + c_1$.  (In this case the 
inequalities $a_1 \geq  0$, $ b_1 \geq  0$, and $c_1 \geq  0$ follow.)  For 
$n \geq  2$, Buch has a precise conjecture for producing a minimal set of 
inequalities: they are the $3n$ inequalities 
$$
	a_1 \geq  a_2 \geq  \ldots  \geq  a_n \geq  0, \,\, 
	 b_1 \geq  b_2 \geq  \ldots  \geq  b_n \geq  0, \,\,
	c_1 \geq  c_2 \geq  \ldots  \geq  c_n \geq  0, 
\tag$**$
$$
and the inequalities ($**_{IJK}$) coming from the following set of triples: 
The triple $(I,J,K)$ must be in $R_r^{2n}$, and the following two 
conditions must be satisfied:
\roster
\item"(i)"  None of $I$, $J$, or $K$ contains a pair of the form 
	$\{t,2n+1-t\}$ for any $1 \leq t \leq n$.
\endroster
This means that there is no cancellation when one writes out the 
corresponding inequality ($**_{IJK}$).  In particular, this requires $r$ to be 
no larger than $n$.  For a subset $I$ of $\{1, \ldots , 2n\}$ of cardinality 
$r$, not containing any pair $\{t, 2n+1-t\}$, set
$$
\overline{I}  = \{ i \in \{1, \ldots ,n\} \,|\, i \in I \text{ or } 2n+1-i \in I\} 
= \{\overline{i}_1 < \overline{i}_2 < \ldots < \overline{i}_r\}, 
$$
and define $\overline{J}$ and $\overline{K}$ similarly from $J$ and $K$.  
The second condition is
\roster
\item"(ii)" 	$\{p \leq r  \, |\,  \overline{k}_p \notin K\} =  \{p 
\leq r  \, |\, \overline{i}_p \notin I\} \,\, \coprod \,\,   
\{p \leq r  \, |\, \overline{j}_p \notin J\}$,
\endroster
where $\coprod$ denotes disjoint union.  This condition 
guarantees that when the 
inequality ($**_{IJK}$) is written out as an inequality with all positive 
coefficients, as we did in the example (24), there are $r$ terms on the left, 
and $2r$ terms on the right, and, for each position $p$ between $1$ and 
$r$, exactly one of the terms $\{a_{\overline{i}_p}, b_{@,\overline{j}_p}, 
c_{@,\overline{k}_p}\}$ is on the left, and the other two are on the right.

\proclaim{Conjecture (Buch)} For $n\geq 2$, the inequalities 
\rom{(}$**$\rom{)} and \rom{(}$**_{IJK}$\rom{)}, for $(I,J,K)$ in 
$R_r^{2n}$, $r \leq n$, satisfying (i) and (ii), define the facets of the cone of 
singular values for $n$ by $n$ matrices.  \endproclaim

Buch has verified this conjecture for $n \leq 4$.  Independently, L. O'Shea 
and R. Sjamaar [OS] have speculated how a complete set of inequalities can 
be obtained by looking at Schubert calculus with $\mathbb{Z}/2\mathbb{Z}$ 
coefficients in real flag varieties.  

\smallpagebreak

In the case of square $n$ by $n$ matrices $A$, $B$, and $C$, one may also 
ask the similar question for products:
\block
	What $a$, $b$, and $c$ can be the singular values of $A$, $B$ and $C$ 
	with $C = A\cdot B$?
\endblock
The answer again is controlled by the triples in $T_r^n$ (or $R_r^n$):
\proclaim{Theorem 16}  A triple $(a,b,c)$ occurs as singular values of $n$ by 
$n$ matrices $A$, $B$, and $C = A\cdot B$ if and only if 
$$
\prod_{k \in K} c_k \leq \prod_{i \in I} a_i \cdot  \prod_{j \in J} b_j 
$$
for all $(I,J,K)$ in $T_r^n$, and all $r < n$. 
\endproclaim
In 1950 Gel'fand and Naimark [GN] proved the special case of this theorem 
when $K = I$ and $J = \{1, \ldots , r\}$.  Many other special cases have 
been found since then, cf. [TT1].

As in the Hermitian case (Theorem 6), these results also extend to singular 
values of compact operators on a Hilbert space. 
 
For the question of which $c_k$, for fixed $k$, can be $k^{\text{th}}$ 
singular values of sums or products of matrices with given singular values, 
the answers are similar to the Hermitian case discussed in Section 1.  For 
these results, see [QS].

\head 6.  First steps toward the proofs \endhead

We now start on the proofs of the theorems, together with some 
generalizations and complements.  In contrast to the first five sections, from 
now on we will prove the propositions as we go along, and include references 
to the original articles.

\subhead 6.1. The Rayleigh trace \endsubhead
An important key to the understanding of the eigenvalue problem was given 
in 1962 by J. Hersch and B. P. Zwahlen, based on then recent work of A. R. 
Amir-Mo\'ez and Wielandt, with ideas going back to Poincar\'e, Weyl, 
Courant, and Fischer.  If $A$ is a Hermitian $n$ by $n$ matrix, with 
eigenvalues $\alpha$, choose corresponding 
orthogonal eigenvectors $v_i$ so $A\,v_i = 
\alpha_i v_i$.  Let $F_k(A)$ be the span of $v_1, \ldots , v_k$, so we have a 
complete flag $F_{\sssize{\bullet}}(A)$ of $\mathbb{C}^n$ (determined 
uniquely by $A$ only if its eigenvalues are distinct).  For any linear subspace 
$L$ of $\mathbb{C}^n$, define the {\bf Rayleigh} \footnote{Lord 
J. W. Rayleigh, a 
British physicist, used the ratios $(Av,v)/(v,v)$ in physical problems involving 
eigenvalues.  For example, the maximum (resp. minimum) of this ratio for all 
nonzero $v$ is the 
biggest (resp. smallest) eigenvalue of $A$.} {\bf trace} $R_A(L)$ by 
$$
	R_A(L) = \sum_{i = 1}^r (A\, u_i, u_i), \tag25
$$
where $u_1, \ldots , u_r$ is an orthonormal basis of $L$.  This is 
independent of the choice of basis; indeed, $R_A(L)$ is the trace of the 
composite $L \to \mathbb{C}^n \to \mathbb{C}^n \to L$, where the first 
map is the inclusion, the second is given by $A$, and the third is orthogonal 
projection.  The key fact, which links the sums occurring in the inequalities 
for the eigenvalues to Schubert varieties, is: 

\proclaim{Proposition 1 [HZ]}  For any subset $I = \{i_1 < \ldots < i_r\}$ of 
$\{1, \ldots , n\}$ of cardinality $r$,
$$
	\sum_{i \in I} \alpha_i = 
	\operatornamewithlimits{Min}_{L \in 
\Omega_I(F_{\sssize{\bullet}}(A))} R_A(L).
$$\endproclaim
The proof is straightforward.  For any $L$ in the Schubert variety 
$\Omega_I(F_{\sssize{\bullet}}(A))$, and for any unit vector $u_1$  in $L\cap 
F_{i_1}(A)$, we have $(A\,u_1,u_1) \geq  \alpha_{i_1}$.  For any unit vector 
$u_2$ perpendicular to $u_1$ in $L\cap F_{i_2}(A)$, we have $(A\,u_2,u_2) 
\geq  \alpha_{i_2}$.  Continuing in this way, one sees that $R_A(L) = 
\sum_{p=1}^r (A\,u_p,u_p) \geq  \sum_{i \in I} \alpha_i$.  This inequality 
will be strict unless these choices can be made with $A\,u_p = 
\alpha_{i_p}u_p$ for $1 \leq p \leq r$.  In particular, taking $L$ spanned 
by $u_p = v_{i_p}$ for $1 \leq p \leq r$, one obtains equality.  

For later use we note the following consequence of this argument.  

\proclaim{Corollary 1}  For $L$ in $\Omega_I(F_{\sssize{\bullet}}(A))$, the 
only way the equality  $\sum_{i \in I} \alpha_i = R_A(L)$  can hold is if $L$ 
is spanned by eigenvectors $u_1, \ldots , u_r$ such that $A(u_p) = 
\alpha_{i_p} u_p$ for $1 \leq p \leq r$.  In this case $A$ maps $L$ \rom{(}and 
therefore also $L^\perp$\rom{)} to itself.  If the eigenvalues $\alpha_1, \ldots , 
\alpha_n$ are distinct, then there is only one $L$ in 
$\Omega_I(F_{\sssize{\bullet}}(A))$ for which $\sum_{i \in I} \alpha_i = 
R_A(L)$.
\endproclaim

The following two corollaries are proved by the same methods, although 
they are not really needed here.

\proclaim{Corollary 2}  Let $F_{\sssize{\bullet}}^\prime(A)$ be the flag 
with $F_k^\prime(A)$ spanned by the last $k$ eigenvectors $v_{n+1-k}, 
\ldots , v_n$, and set $I^\prime = \{i \, | \, n+1-i \in I\}$. Then 
$$
	\sum_{i \in I} \alpha_i = 
	\operatornamewithlimits{Max}_
	{L \in \Omega_{I^\prime}(F_{\sssize{\bullet}}^\prime(A))} \, R_A(L).
$$\endproclaim
To prove this, construct a sequence of orthonormal vectors $u_p$ in $L\cap 
F_{i_p^\prime}^\prime(A)$.  One has $(A\,u_p,u_p) \leq \alpha_{i_{r+1-p}}$ 
for $1 
\leq p \leq r$.  Taking $L$ spanned by $v_{i_1}, \ldots , v_{i_r}$ 
gives equality.

\proclaim{Corollary 3}  Let $A_L : L \to L$ denote the composite $L \to 
\mathbb{C}^n \to \mathbb{C}^n \to L$ defined before the proposition, and 
let $\alpha_1^L \geq  \ldots  \geq  \alpha_r^L$ be its eigenvectors.  Then 
$$
\alpha_{n-r+k} \leq \alpha_k^L \leq \alpha_k \quad 
 \text{ for } \quad  1 \leq k \leq r.
$$\endproclaim
For this, let $w_1, \ldots , w_r$ be an orthonormal basis of $L$ with 
$A_L(w_i) = \alpha_i^Lw_i$, $1 \leq i \leq r$.  The intersection 
$$
   	\langle w_1, \ldots , w_k \rangle \cap \langle v_1, \ldots , v_{k-1} 
	\rangle ^\perp = \langle w_1, \ldots , w_k \rangle \cap \langle v_k, 
	\ldots , v_n\rangle 
$$
is not empty, and if $v$ is a unit vector in it, then $\alpha_k^L \leq (A\,v,v) 
\leq \alpha_k$.  Taking a unit vector in $\langle w_k, \ldots , w_r \rangle 
\cap \langle v_1, \ldots , v_{n-r+k} \rangle$ gives the other inequality.

When $L$ is the subspace spanned by some of the basic vectors of 
$\mathbb{C}^n$, the matrix for $A_L$ is a principal minor of $A$, and the 
result of Corollary 3 is sometimes called the {\em inclusion principle}.

With a little Schubert calculus, Proposition 1 leads quickly to a proof of the 
following proposition.  

\proclaim{Proposition 2  ([J], [HR], [K], [T])}  Let $A$, $B$, and $C$ be 
Hermitian $n$ by $n$ matrices, with $C = A + B$, and eigenvalues 
$\alpha$, $\beta$, and $\gamma$.  Then for any $r < n$ and any $(I,J,K)$ 
in $S_r^n$, the inequality \rom{(}$*_{IJK}$\rom{)} is valid. \endproclaim
To prove this, note first that the Rayleigh trace is linear in the matrices, so 
$$
	R_{-A}(L) + R_{-B}(L) + R_C(L) = R_0(L) = 0. \tag26
$$
Note also that if $\alpha_1 \geq  \ldots  \geq  \alpha_n$ are the 
eigenvalues of $A$, then 
$\alpha_1^\prime = -\alpha_n \geq  \ldots  \geq  \alpha_n^\prime = -
\alpha_1$ are the eigenvalues of $-A$.  For a subset $I$ of $\{1, \ldots 
,n\}$, let $I^\prime = \{i \,|\, n+1-i \in I\}$. This gives the equations 
$$
 \sum_{i \in I^\prime} \alpha_i^\prime = - \sum_{i \in I} \alpha_i . \tag27
$$
It follows from the definition of Schubert classes (see (21)) that 
$$
	\sigma_{\lambda(I)} = \omega_{I^\prime}, \tag28
$$
where the notation is that of Section 4, but working now in the 
Grassmannian $Gr(r,\mathbb{C}^n)$.  In addition,
$$
\sigma_{\lambda(I)} \,\, \text{ and } \,\, \sigma_{\lambda(I^\prime)} 
\,\, \text{ are 
dual}. \tag29
$$
This means that $\sigma_{\lambda(I)} \cdot \sigma_{\lambda(J^\prime)} = 
\sigma_\rho$ if and only if $I = J$, where $\rho = (n-r)^r$, so 
$\sigma_\rho$ is the class of a point.  
The following are therefore equivalent:  
\roster
\item"(i)"  $\sigma_{\lambda(K)}$ occurs in $\sigma_{\lambda(I)} \cdot 
		\sigma_{\lambda(J)}$, i.e., $(I,J,K) \in S_r^n$;  
\item"(ii)" $\sigma_{\lambda(I)} \cdot \sigma_{\lambda(J)} \cdot 
		\sigma_{\lambda(K^\prime)} \ne 0$;  
\item"(iii)" $\omega_{I^\prime} \cdot \omega_{J^\prime} \cdot \omega_K 
		\ne 0$.  
\endroster
This last condition implies that the intersection  
$\Omega_{I^\prime}(F_{\sssize{\bullet}}(-A)) \cap  
\Omega_{J^\prime}(F_{\sssize{\bullet}}(-B)) \cap  
\Omega_K(F_{\sssize{\bullet}}(C))$ is not empty.  Taking $L$ in this 
intersection, this gives
$$
	-  \sum_{i \in I} \alpha_i  -  \sum_{j \in J} \beta_j + \sum_{k \in K} 
	\gamma_k \, \leq \,  R_{-A}(L) + R_{-B}(L) + R_C(L) = 0,
$$
and this proves Proposition 2.

\remark{Remark}  Using the equation $A + B + (-C) = 0$, one deduces 
by a similar argument that $\sum_{k \in K'} \gamma_k \leq \sum_{i \in 
I'} \alpha_i + \sum_{j \in J'} \beta_j$,  where 
$I' = \{i \,|\, n+1-i \notin I\}$.  In fact, however, a triple $(I,J,K)$ 
is in $S_r^n$ if and only if $(I',J',K')$ is 
in $S_{n-r}^n$.  This follows from the isomorphism between 
$Gr(r,\mathbb{C}^n)$ and $Gr(n-r,\mathbb{C}^n)$ that takes a subspace 
$L$ of $\mathbb{C}^n$ to the kernel of the map $\mathbb{C}^n = 
(\mathbb{C}^n)^* \to L^*$.  On cohomology this map takes 
$\sigma_\lambda$  to $\sigma_{\lambda'}$ , where 
$\lambda'$  is the conjugate of $\lambda$.  Note that 
$\lambda(I')$ is the transpose of $\lambda$ (see 
Lemma 4 in Section 10).  \endremark

In 1962, and for some time after that, people working on this problem were 
not familiar with Schubert calculus, and it was only in 1979, in an 
unpublished thesis of S. Johnson [J], directed by Thompson, that the 
connection with Schubert calculus was made.   Proposition 2 was later 
rediscovered by U. Helmke and J. Rosenthal [HR], Totaro [T] (in a 
different context), and Klyachko [Kl1].  Many papers proved special cases of 
the inequalities ($*_{IJK}$) by showing by hand --- using only linear algebra 
--- that the intersection of Schubert varieties appearing in the above proof 
must be nonempty.

\subhead 6.2. Relations among the four main problems \endsubhead
The fact that the representation theory problem is equivalent to the 
Schubert calculus problem has been known at least since 1947, when L. 
Lesieur [Le] proved that the Littlewood-Richardson coefficients 
$c_{\alpha\,\beta}^{\,\,\gamma}$  are the same as the coefficients 
$d_{\alpha\,\beta}^{\,\,\gamma}$ that describe the multiplication in the 
cohomology of the Grassmannian.  In spite of all the known relations 
between Schubert varieties and representation theory, however, the proof is 
not very direct.  It proceeds by showing that both are controlled by the same 
algebra of symmetric functions.  The character of the representation 
$V(\alpha)$ of ${GL}_n(\mathbb{C})$ is the Schur polynomial 
$s_\alpha(x_1, \ldots , x_n)$.  Corresponding to the decomposition of a 
tensor product one has the identity $s_\alpha(x) \cdot s_\beta(x) = \sum  
c_{\alpha\,\beta}^{\,\,\gamma} \, s_\gamma(x)$.  These Schur polynomials 
form a basis for the ring $\Lambda_n$ of symmetric polynomials in the 
variables $x_1, \ldots , x_n$ (cf. [Mac], [Fu2]).  
On the other hand, there is a 
surjection of the ring $\Lambda_n$ onto the cohomology of the 
Grassmannian $Gr(n,\mathbb{C}^m)$, that takes $s_\alpha(x)$ to the 
class $\sigma_\alpha$.  The classical formulas of Pieri and Giambelli imply 
that this is a ring homomorphism, which proves the following proposition 
(see e.g. [Fu2] for this story):

\proclaim{Proposition 3}  The coefficients for multiplying Schubert classes 
are the same as those for decomposing tensor products, i.e., 
$d_{\alpha\,\beta}^{\,\,\gamma} = c_{\alpha\,\beta}^{\,\,\gamma}$.
\endproclaim

Closer relations are known between the cohomology of the Grassmann 
variety $G/P$ and representations of $G$ stemming from work of 
Ehresmann and Kostant (see [Ko],\S8), but one would like a more direct 
explanation for the equality of Proposition 3.

Based on the work of Hall and Green, in 1968 Klein [K] completed the 
proof of the following 
proposition:

\proclaim{Proposition 4}  If $R$ is a discrete valuation ring, and 
$\mathcal{C}$ is an $R$-module with invariant factors $\gamma$, then there 
is a submodule $\mathcal{B}$ of $\mathcal{C}$ with invariant factors 
$\beta$  such that $\mathcal{A} = \mathcal{C}/\mathcal{B}$ has invariant 
factors $\alpha$  if and only if the Littlewood-Richardson coefficient 
$c_{\alpha\,\beta}^{\,\,\gamma}$ is positive. \endproclaim

Although this establishes a relation between invariant factors and 
representation theory, the proof does nothing of the sort.  
Rather, it provides a link between such submodules and the 
Littlewood-Richardson algorithm.  In fact, if $\mathcal{B}$ is 
a submodule of $\mathcal{C}$, the invariant factors $\gamma(k)$ of 
$\mathcal{C}/\pi^k\mathcal{B}$ form a chain of partitions 
$\alpha = \gamma(0) \subset \gamma(1) \subset \ldots \subset \gamma(t) 
= \gamma$.  If a $k$ is placed in each box of $\gamma(k) \smallsetminus 
\gamma(k-1)$, Green showed that one obtains an array whose conjugate 
--- obtained by interchanging rows and columns --- satisfies the conditions 
(i) -- (iv) of Section 2 for the conjugate partitions $\alpha'$, $\beta'$, 
$\gamma'$.  Conversely, Klein showed that any such chain of partitions 
can be realized by some submodule $\mathcal{B}$ of $\mathcal{C}$.  To 
complete the proof, one needs the fact that $c_{\alpha'\,\beta'}^{\,\,\gamma'} 
= c_{\alpha\,\beta}^{\,\,\gamma}$ (see Lemma 4 in Section 10).

There are polynomials $g_{\alpha\,\beta}^{\,\,\gamma}(T)$, called Hall 
polynomials, whose top coefficient is $c_{\alpha\,\beta}^{\,\,\gamma}$, 
such that whenever the residue field of $R$ is finite of cardinality $q$, 
then $g_{\alpha\,\beta}^{\,\,\gamma}(q)$ is the number of submodules 
$\mathcal{B}$ of $\mathcal{C}$ with invariant factors $\beta$, such that 
$\mathcal{C}/\mathcal{B}$ has invariant factors $\alpha$.
For a proof, with 
references, applications, and history, see Chapter II of [Mac].

Thompson and his coworkers, who worked on both problems, noticed that 
the answers to the Hermitian eigenvalue problem and the invariant factor 
problem appeared to be similar to each other, and were therefore related to 
Littlewood-Richardson coefficients by means of Proposition 4. 

Although this was only realized in 1997 by A. P. Santana, J. F. Queir\'o, and 
E.~Marques de S\'a [SQS], the following relation can be deduced from a 1982 
theorem of G.~J.~Heckman [He] about moment mappings and coadjoint 
orbits:

\proclaim{Proposition 5}  If $\alpha$, $\beta$, and $\gamma$  are 
partitions of lengths at most $n$ such that $c_{\alpha\,\beta}^{\,\,\gamma} 
\ne 0$, then there are Hermitian $n$ by $n$ matrices $A$, $B$, and $C$ 
with eigenvalues $\alpha$, $\beta$, and $\gamma$  and $C = A + B$.
\endproclaim

In fact, we will not need this proposition; it will follow from the other 
results proved here. 
We will discuss the coadjoint orbit and moment map approach in Section 8. 
Note, however, that with these five propositions, 
relations are established among all four of our subjects.  For partitions  
$\alpha$, $\beta$, and $\gamma$  of lengths at most $n$, consider the 
following conditions: 
\roster
\item"(i)" $(\alpha,\beta,\gamma)$ occur as 
eigenvalues of $n$ by $n$ Hermitian matrices $A$, 
$B$, and $C = A + B$; 
\item"(ii)" $(\alpha,\beta,\gamma)$ occur as 
invariant factors (over one or every discrete 
valuation ring) of $n$ by $n$  matrices $A$, $B$, and $C = A\cdot B$; 
\item"(iii)" the representation $V(\gamma)$ of ${GL}_n(\mathbb{C})$ 
occurs in $V(\alpha)\otimes V(\beta)$; 
\item"(iv)" $\sigma_\gamma$  occurs in $\sigma_\alpha \cdot 
\sigma_\beta$  in the cohomology of $Gr(n,\mathbb{C}^m)$ for any large 
$m$.  
\endroster
Propositions 3 and 4 say that (ii), (iii), and (iv) are equivalent, and 
Proposition 5 says that they imply (i).  See [SQS] for more about these 
relations. 

Note that the Littlewood-Richardson coefficients measure a multiplicity in 
each of the problems (ii), (iii), and (iv).  It would be interesting to find a 
similar interpretation for the eigenvalue problem (i). 
\footnote{Knutson reports that there is an 
interpretation of a suitable volume of a space of Hermitian matrices 
$(A,B,C)$ with $C=A+B$ and eigenvalues $\alpha$, $\beta$, and $\gamma$, in 
terms of Littlewood-Richardson coefficients, but asymptotically, using 
triples $(N\alpha,N\beta,N\gamma)$ for large $N$. See also [DRW], \S3.}

Proposition 2 established a link in one direction between the eigenvalue 
problem for $n$ by $n$ matrices and Schubert calculus in smaller 
Grassmann varieties $Gr(r,\mathbb{C}^n)$.  The next main step is to find a 
converse to this proposition.  It should be pointed out at this point 
that it is 
not obvious that {\em any} of the four problems have solutions that can be 
described by any inequalities of the form ($*_{IJK}$).  The fact that the 
eigenvalue problem can be described by some linear inequalities also follows from 
convexity properties of moment maps, which we will discuss in Section 8.

\head 7.  Filtered vector spaces, geometric invariant theory, and stability \endhead  

We start by describing a theorem of Totaro [T], which, although written for 
quite another purpose, quickly yields the converse to Proposition 2.  I thank 
R. Lazarsfeld for pointing me to this work of Totaro.  Let $V$ be an 
$n$-dimensional 
complex vector space.  By a {\bf filtration} 
$V^{\sssize{\bullet}}$  of $V$ we shall mean a weakly decreasing sequence 
of subspaces
$$
	V = V^0 \supset V^1 \supset V^2 \supset \ldots \supset V^p \supset 
	V^{p+1} \supset \ldots 
$$
with the assumption that $V^q = 0$ for sufficiently large $q$.  Any subspace 
$L$ of $V$ gets an induced filtration, by setting $L^p = L \cap V^p$.

By an $m$-{\bf filtration} of $V$ is meant a collection 
$V^{\sssize{\bullet}}(s)$ of filtrations of $V$, for $1 \leq s \leq m$.  For any 
nonzero subspace $L$ of $V$, define its {\bf slope} $\mu(L)$ with respect 
to this $m$-filtration by 
$$ 
	\mu(L) = \frac1{\operatorname{dim}(L)} \sum_{s=1}^m \sum_{ p 
\geq 1} \operatorname{dim}(L\cap V^p(s)). \tag30
$$
Call the $m$-filtration {\bf semistable} (resp. {\bf stable}) if $\mu(L) \leq 
\mu(V)$ (resp. $\mu(L) < \mu(V)$) for all subspaces $L$ of $V$, $0 \ne L 
\ne V$.  There is an obvious notion of a direct sum of $m$-filtered vector 
spaces, taking the $p^{\text{th}}$ subspace of the sum to be the sum of the 
$p^{\text{th}}$ subspaces of the factors.  An $m$-filtration is called {\bf 
polystable} if it 
is a direct sum of a finite number of stable $m$-filtrations of 
the same slope.

\proclaim{Proposition 6 [T]}  An $m$-filtration on $V$ is polystable if and 
only if there is a Hermitian metric on $V$ such that the sum of the 
orthogonal projections from $V$ onto the spaces $V^p(s)$, for $p \geq1$ 
and $1 \leq s \leq m$, is the scalar operator $\mu(V)$.  \endproclaim
The proof of this uses some geometric invariant theory.  Totaro shows that a 
polystable $m$-filtration corresponds to a polystable point $x$ in a product 
$X$ of $m$ partial flag varieties, which has a canonical embedding in a 
projective space $\mathbb{P}(E)$, where $E = 
\otimes_{s=1}^m\otimes_{p\geq1}\bigwedge^{\operatorname{dim} 
V^p(s)} V$.  It is a result of G. Kempf and L. Ness from Geometric Invariant 
Theory, that if $\widetilde{x}$ is a representative of $x$ in $E$, 
polystability is equivalent to the corresponding orbit ${SL}(V) \cdot 
\widetilde{x}$ being closed in $E$.  One then chooses a Hermitian metric 
on $V$ (which induces a metric on $E$) to minimize the distance from 
$\widetilde{x}$ to the origin.  The fact that this is a critical point implies, 
by a calculation, the asserted fact that the sum of the projections is a 
scalar.

This proposition was used by Totaro to give a simpler proof of a theorem of 
G. Faltings, that the tensor product of semistable filtrations is semistable.  
Essentially the same argument was found by Klyachko [Kl1] for the eigenvalue 
problem.  See Section 8 and [Fu3] for more details.  

If $\alpha  = (\alpha_1 \geq  \alpha_2 \geq  \ldots  \geq  \alpha_n \geq  
0)$ is a partition of length at most $n$, and if $F_{\sssize{\bullet}} : 0 = 
F_0 
\subset F_1 \subset \ldots \subset F_n = V$ is a complete flag in $V$, 
define a filtration $V^{\sssize{\bullet}}$ of $V$ by setting $V^p = 
F_{{\alpha'}_p}$, where $\alpha'$  is the conjugate 
partition to $\alpha$, i.e., ${\alpha'}_p$ is the cardinality of the 
set $\{i \in \{1, \ldots , n\} 
\,|\, \alpha_i \geq  p\}$.  Any filtration arises 
in this way, for some unique $\alpha$  and some flag $F_{\sssize{\bullet}}$, 
although the flag is not unique unless $\alpha$  consists of $n$ distinct 
integers.  The following lemma is immediate from this definition.

\proclaim{Lemma 1}  Choose a Hermitian metric on $V$, and let $A$ be the 
sum of the orthogonal projections of $V$ on $V^p$, for $p \geq 1$.  Then 
$A$ is a Hermitian operator whose eigenvalues are $\alpha_1, \ldots , 
\alpha_n$. \endproclaim

For any subset $I = \{i_1 < \ldots < i_r\}$ of $\{1, \ldots , n\}$ of 
cardinality $r$, we have the {\bf Schubert cell} $\Omega_I^\circ 
(F_{\sssize{\bullet}})$, defined by
$$
\Omega_I^\circ (F_{\sssize{\bullet}}) = \{ L \in Gr(r,V) \,\,| \,\,
\operatorname{dim}(L\cap F_j) = k \, \text{ for }\, i_k \leq j < i_{k+1}, 
\, 0 \leq k \leq r\},
$$
where 
$i_0$ is defined to be $0$ and 
$i_{r+1}$ is defined to be $n+1$; 
in other words, $I$ lists those $i$ for which $L\cap F_i \neq L\cap F_{i-1}$.  
This cell is a dense open subset of the 
Schubert variety $\Omega_I(F_{\sssize{\bullet}})$.  Any subspace belongs to 
a unique Schubert cell, when the flag $F_{\sssize{\bullet}}$ is fixed.  From 
this definition it follows that, for any $\alpha_1, \ldots , \alpha_n$, 
$$
\sum_{i \in I} \alpha_i = 
\sum_{i=1}^n\alpha_i \operatorname{dim}(L\cap F_i/L\cap F_{i-1}) 
= \sum_{i=1}^n(\alpha_i - \alpha_{i+1}) \operatorname{dim}(L\cap F_i), 
\tag31
$$
where $\alpha_{n+1}$ is defined to be $0$.  By a simple 
counting argument, this implies: 

\proclaim{Lemma 2}  Let $L$ be a subspace of $V$ that is in the Schubert 
cell $\Omega_I^\circ (F_{\sssize{\bullet}})$.  Then $\sum_{i \in I} 
\alpha_i = \sum_{p \geq 1} \operatorname{dim}(L\cap V^p)$. 
\endproclaim

We define $S_r^n(m)$ to be the set of $m$-tuples $\mathscr{I} = (I(1), 
\ldots , I(m))$ such that $\prod_{s=1}^m \omega_{I(s)} \ne 0$.  Here, as in 
Section 4, $\omega_I$ denotes the class of a Schubert variety 
$\Omega_I(F_{\sssize{\bullet}})$. 

Consider now $m$-tuples $\alpha(1), \ldots , \alpha(m)$, with each 
$\alpha(s)$ a nonincreasing sequence of real numbers of length at most 
$n$, written $\alpha(s) = (\alpha_1(s), \ldots , \alpha_n(s))$.  For any 
$m$-tuple $\mathscr{I} = (I(1), \ldots , I(m))$ of subsets of $\{1, \ldots , 
n\}$ of cardinality $r$ we have the corresponding inequality
$$
\frac1{r} \sum_{s=1}^m \sum_{i \in I(s)} \alpha_i(s)  \,\, \leq \, \,\, 
\frac1{n} \sum_{s=1}^m \sum_{i=1}^n \alpha_i(s).
\tag$*_{\mathscr{I}}$
$$

\proclaim{Lemma 3}  Let $\alpha(s)$ be a partition of length at most $n$, 
for $1 \leq s \leq m$.  Let $F_{\sssize{\bullet}}(s)$ be general flags in $V$, 
$1 \leq s \leq m$, and let $V^{\sssize{\bullet}}(s)$ be the filtration defined 
by the flag $F_{\sssize{\bullet}}(s)$ and the partition $\alpha(s)$, for $1 
\leq s \leq m$.  
This $m$-filtration is semistable (resp. stable) if and only if 
the inequality  \rom{(}$*_{\mathscr{I}}$\rom{)} holds (resp. with strict 
inequality) for every  $\mathscr{I}$ in $S_r^n(m)$ and all $r \leq n$. 
\endproclaim
This follows from Lemma 2 and the fact that, when the 
$F_{\sssize{\bullet}}(s)$ are general flags, an intersection of Schubert cells 
$\bigcap_{s=1}^m \Omega_{I(s)}^\circ (F_{\sssize{\bullet}}(s))$ is 
nonempty if and only if the corresponding class $\prod_{s=1}^m 
\omega_{I(s)}$ is not zero.

The following, in the Hermitian case, is one of the main results of 
Klyachko's paper [Kl1].  
From now on, we will state the results for an arbitrary 
number of factors, rather than the three that were featured in the first 
section; and we will put them all on the same side of the equation.  This 
allows a simpler and more natural expression of the results and proofs, by 
avoiding the kind of manipulations that we saw in the proof of Proposition 2.  

\proclaim{Proposition 7}  Let $\alpha(1), \ldots , \alpha(m)$ be weakly 
decreasing sequences of $n$ real numbers.  There are Hermitian 
$n$ by $n$ matrices $A(1), \ldots , A(m)$, with $\alpha(s)$ the 
eigenvalues of $A(s)$, and $\sum_{s=1}^m A(s)$ a scalar matrix, if and only 
if \rom{(}$*_{\mathscr{I}}$\rom{)} holds for every $\mathscr{I}$ in 
$S_r^n(m)$, for all $r < n$. \endproclaim

The implication $\Rightarrow$ is proved as in Proposition 2, but simpler:  If 
$\sum_sA(s) = c$, and $\prod_s \omega_{I(s)} \ne 0$, there is 
some subspace $L$ in $\bigcap_s \Omega_{I(s)}^\circ 
(F_{\sssize{\bullet}}(s))$.  By Proposition 1, 
$$
\sum_{s=1}^m \sum_{i \in I(s)} \alpha_i(s) \leq \sum_{s=1}^m  
R_{A(s)}(L) = R_c(L) = c \cdot r, 
$$
and $c = \frac1{n}\operatorname{Trace}(\sum_s A(s)) = \frac1{n} \sum_s 
\sum_{i=1}^n \alpha_i(s)$.  This proves that ($*_{\mathscr{I}}$) holds.  

For the converse, first one checks that the region defined by the inequalities 
($*_\mathscr{I}$) and the inequalities that make each of $\alpha(s)$ 
weakly decreasing has a nonempty interior ([Fu3], Lemma 2).  By a 
continuity argument, using the compactness of the unitary group $U(n)$, it 
suffices to prove the existence of the Hermitian matrices when each 
$\alpha(s)$ consists of $n$ distinct rational numbers, and all of the 
inequalities ($*_\mathscr{I}$) are strict; therefore (since multiplying 
matrices by scalars or adding scalar matrices doesn't change the situation), 
we may assume that each $\alpha(s)$ consists of $n$ distinct nonnegative 
integers.  Take general flags $F_{\sssize{\bullet}}(s)$, $1 \leq s \leq m$, 
and use the partitions $\alpha(s)$ to construct filtrations 
$V^{\sssize{\bullet}}(s)$.  We have seen that this $m$-filtration is stable.  By 
Proposition 6, there is a Hermitian metric on $V$ so that the sum of the 
projections on the spaces $V^p(s)$, $p \geq 1$, $1 \leq s \leq m$, is a 
scalar.  The conclusion follows from Lemma 1. 

Special cases of the following proposition were stated by V. B. Lidskii [L1] 
and proved by Horn [H2]. In fact, Horn used only calculus to show that, 
assuming that the 
eigenvalues are distinct, if one is on a boundary of the region of 
possible eigenvalues for real symmetric matrices, there must be an 
invariant subspace. This was an early step toward his conjecture. 

\proclaim{Proposition 8}  Let $A(s)$ be Hermitian $n$ by $n$ matrices, 
with eigenvalues $\alpha(s)$, for $1 \leq s \leq m$, whose sum $\sum_s 
A(s)$ is a scalar.  Suppose that there is some $\mathscr{I}$ in $S_r^n(m)$ 
such that \rom{(}$*_\mathscr{I}$\rom{)} holds with equality.  Then there 
is an 
$r$-dimensional subspace $L$ of $\mathbb{C}^n$ that is mapped to itself by 
each $A(s)$. \endproclaim
To prove this, take a flag $F_{\sssize{\bullet}}(s) = 
F_{\sssize{\bullet}}(A(s))$ corresponding to $A(s)$ as in Proposition 1.  Let 
$c = \frac1{n} \sum_s \sum_i \alpha_i(s)$ be the scalar that the matrices 
$A(s)$ sum to.  Since $\mathscr{I}$ is in $S_r^n(m)$ there must be an 
$r$-dimensional subspace $L$ of $\mathbb{C}^n$ in 
$\bigcap_s\Omega_{I(s)}(F_{\sssize{\bullet}}(s))$.  By Proposition 1,
$$
	\frac1{r} \sum_s \sum_{i \in I(s)} \alpha_i(s) \leq \frac1{r} 
	\sum_s R_{A(s)}(L) = \frac1{r} R_{\sum_sA(s)}(L) = c.
$$
Since ($*_\mathscr{I}$) holds with equality, we must have $\sum_{i \in 
I(s)} \alpha_i(s) = R_{A(s)}(L)$ for every $s$.  By Corollary 1 to Proposition 
1, each $A(s)$ must map $L$ to itself.

Since each $A(s)$ is Hermitian, it follows that it also maps $L^\perp$ to 
itself.  The sum of the restrictions of $A(s)$ to $L$ (or $L^\perp$) is the 
same scalar $c$, so the process can be repeated for these: either all the 
inequalities for them are strict, or they can be further decomposed.  (This is 
essentially Totaro's proof [T] of Proposition 6 above.) Proposition 8 can also 
be proved from properties of moment maps, as in [Kn].

\proclaim{Corollary}  Suppose in addition that each $A(s)$ is a real 
symmetric matrix, and $\prod_s \omega_{I(s)}$ is an odd multiple of 
$\omega_{\{1,\ldots,r\}}$, where $r$ is the cardinality of each $I(s)$.  
Then there is an $r$-dimensional subspace $L$ of 
$\mathbb{R}^n$ that is preserved by each $A(s)$. \endproclaim
Indeed, the fact that $\prod_s \omega_{I(s)}$ is an odd multiple of 
$\omega_{\{1,\ldots,r\}}$ guarantees that the intersection $\bigcap_s 
\Omega_{I(s)}(F_{\sssize{\bullet}}(s))$ of real Schubert varieties must 
contain a real point.  
This is obvious if the flags are in general position, since 
the complex points occur in pairs.  It is also true for arbitrary flags, by the 
results of [Fu1], Chapter 13. 

This explains Example 1 in Section 1, which presents  
real symmetric matrices for which 
an equality is satisfied but for which 
there is no real subspace preserved by the 
matrices.  In the present terminology, each $I(s) = \{2,4,6\}$, and 
$\prod_s \omega_{I(s)} = 2\omega_{\{1,\ldots,r\}}$; the corresponding 
intersection of three real Schubert varieties has two complex points, but no 
real points. See Sottile's article [So] for general results about real 
Schubert calculus.   

If each $A(s)$ is real symmetric, if at least one has $n$ distinct 
eigenvalues, and if some ($*_\mathscr{I}$) holds with equality, then there is a 
real subspace of dimension $r$ preserved by each of the $A(s)$.  Indeed, 
the fact that the eigenvalues of 
at least one $A(s)$ are distinct guarantees (by 
Proposition 10 below) that $\prod_{s=1}^m 
\omega_{I(s)} = \omega_{\{1,\ldots,r\}}$, 
so the corollary applies.   

It is not hard to see that, in 
the list of inequalities ($*_\mathscr{I}$), those 
for which the sum of the codimensions of the $\omega_{I(s)}$ is less than 
$r(n-r)$ follow from those where the sum of the codimensions is equal to 
$r(n-r)$.  The latter are the classes for which $\prod_s \omega_{I(s)}$ is a 
nonzero multiple of the class $\omega_{\{1,\ldots,r\}}$ of a point.
This is a consequence of Pieri's formula for multiplying a 
Schubert class by the class $\sigma_{(1)}$ of a hyperplane.  This also follows 
from the stronger results below.  

C. Woodward was the first to discover that this reduced list of inequalities 
($*_\mathscr{I}$) is still redundant, by finding the example 
(12).  These examples involve Schubert classes whose intersection number 
is greater than $1$.  P. Belkale conjectured and proved 
that all such inequalities 
must be redundant.  To state his theorem, let  
$R_r^n(m)$ be the set of 
$m$-tuples
$\mathscr{I}$ of subsets of $\{1, \ldots ,n\}$ of cardinality $r$ such 
that $\prod_s \omega_{I(s)} = \omega_{\{1,\ldots,r\}}$. 

\proclaim{Proposition 9}  On the space of $m$-tuples $\alpha(1), \ldots , 
\alpha(m)$ of weakly decreasing sequences of $n$ real numbers, the 
inequalities \rom{(}$*_\mathscr{I}$\rom{)} for all $\mathscr{I}$ in 
$S_r^n(m)$, and all $r < n$, are implied by those in 
$R_r^n(m)$ for all $r < n$. 
\endproclaim

In fact, Knutson and Woodward proved that the extra inequalities are 
redundant in a 
particularly strong way, as we had observed in examples:

\proclaim{Proposition 10}  Let $\alpha(1), \ldots , \alpha(m)$ be weakly 
decreasing sequences of $n$ real numbers, satisfying the inequalities 
\rom{(}$*_\mathscr{I}$\rom{)} for $\mathscr{I}$ in $R_r^n(m)$, all $r < 
n$.  Suppose, for some $r < n$ and some $\mathscr{I}$ in $S_r^n(m) 
\smallsetminus R_r^n(m)$, we have the reverse inequality
$$
  	\frac1{r} \sum_{s=1}^m \sum_{i \in I(s)} \alpha_i(s) \, \geq \,  
	\frac1{n} \sum_{s=1}^m \sum_{i=1}^n \alpha_i(s).
$$
Then this inequality is actually an equality.  Moreover, for every $s$ between 
$1$ and $m$, at least one of the inequalities $\alpha_1(s) \geq  \alpha_2(s) 
\geq  \ldots  \geq  \alpha_n(s)$ must be an equality.  \endproclaim

The proofs of Belkale and Woodward follow an idea familiar in the study of 
stability of vector bundles: a maximal destabilizing subbundle with maximal 
slope must be unique.  As before, take general complete flags 
$F_{\sssize{\bullet}}(s)$ in $V = \mathbb{C}^n$.  For any $r$-dimensional 
subspace $L$ of $V$, there are unique subsets $I(s)$ of $\{1, \ldots , n\}$ 
of cardinality $r$ such that $L$ is in the intersection $\bigcap_{s=1}^m 
\Omega_{I(s)}^\circ (F_{\sssize{\bullet}}(s))$.  Let $\mathscr{I}(L) = (I(1), 
\ldots , I(m))$ be this 
$m$-tuple of subsets.  Note that, since the flags are general, this 
intersection of Schubert cells is nonempty if and only if the corresponding 
product $\prod_{s=1}^m\omega_{I(s)}$ of Schubert classes is not zero.  
Define the {\bf slope} of $L$, $\mu(L)$, by the formula
$$
\mu(L) = \frac1{r} \sum_{s=1}^m \sum_{i \in I(s)} \alpha_i(s). \tag32
$$
Suppose $\mathscr{I}$ is in $S_r^n(m) \smallsetminus R_r^n(m)$ and 
$\frac1{r} \sum_{s=1}^m \sum_{i \in I(s)} \alpha_i(s)  >  \frac1{n} 
\sum_{s=1}^m \sum_{i=1}^n \alpha_i(s)$.  There must be some $L$ in the 
intersection $\bigcap_{s=1}^m \Omega_{I(s)}^\circ 
(F_{\sssize{\bullet}}(s))$, and $\mu(L) > \mu(V)$ for any such $L$.  Let 
$\mu$ be the maximum of all slopes of all subspaces $L$, and take $r$ 
maximal such that some $L$ of dimension $r$ has slope $\mu$.  Choose a 
subspace $L$ of dimension $r$ with $\mu(L) = \mu$, and set $\mathscr{I} 
= \mathscr{I}(L)$.  Since $\mu(L) > \mu(V)$, $\mathscr{I}$ is not in 
$R_r^n(m)$, so the intersection $\bigcap_{s=1}^m \Omega_{I(s)}^\circ 
(F_{\sssize{\bullet}}(s))$ has cardinality greater than $1$.  Let  $L^\prime 
\ne L$ be another space in the same intersection.  Consider the sum $L + 
L^\prime$, of some dimension $r+t > r$, and set $\mathscr{K} = 
\mathscr{I}(L+L^\prime)$.  The intersection $L\cap L^\prime$ has 
dimension $r-t$; set $\mathscr{J} = \mathscr{I}(L\cap L^\prime)$.  

The linear algebra fact that, for all $i$ and $s$,
$$
   	\operatorname{dim}((L+L^\prime)\cap F_i(s)) \geq  
	\operatorname{dim}(L\cap F_i(s)) + 
	\operatorname{dim}(L^\prime\cap F_i(s)) - 
	\operatorname{dim}((L\cap L^\prime)\cap F_i(s)),
$$
together with equation (31), implies that 
$$
 	\sum_{k \in K(s)} \alpha_k(s)\, \geq  \sum_{i \in I(s)} \alpha_i(s) + 
	\sum_{i \in I(s)} \alpha_i(s) - \sum_{j \in J(s)} \alpha_j(s).
$$
Since $r \cdot \mu(L) = \sum_s \sum_{i \in I(s)} \alpha_i(s)$, and 
similarly for the others, this implies the inequality 
$$
	(r+t) \cdot \mu(L+L^\prime) \,\geq \,  r \cdot \mu(L) + r\cdot \mu(L) - 
	(r-t) \cdot \mu(L\cap L^\prime).
$$
Since $\mu(L) = \mu$ and $\mu(L\cap L^\prime) \leq \mu$, the right 
side is at least $2r\cdot \mu - (r-t) \cdot \mu = (r+t)\cdot \mu$.  Hence 
$\mu(L+L^\prime) \geq  \mu$, which violates the maximality of $r$.  
See [Gr] for more about maximal destabilizing subspaces.

For the last assertion in Proposition 10, suppose that all the inequalities 
($*_\mathscr{I}$) hold for $\mathscr{I}$ in $S_r^n(m)$, but that there is 
some $\mathscr{I}$ in $S_r^n(m) \smallsetminus R_r^n(m)$ for which 
equality holds.  Suppose for some $s$ that the eigenvalues $\alpha_1(s), 
\ldots , \alpha_n(s)$ are all distinct.  We know that, for an appropriate 
metric on $V = \mathbb{C}^n$, there are Hermitian operators $A(s)$ whose 
sum is scalar and so that each of the general flags $F_{\sssize{\bullet}}(s)$ is 
a flag $F_{\sssize{\bullet}}(A(s))$ corresponding to $A(s)$.  Since 
$\mathscr{I}$ is not in $R_r^n(m)$, there are at least two $r$-dimensional 
linear subspaces in $\bigcap_{s=1}^m 
\Omega_{I(s)}(F_{\sssize{\bullet}}(s))$.  
But now for any $L$ in $\bigcap_{s=1}^m 
\Omega_{I(s)}(F_{\sssize{\bullet}}(s))$, by Proposition 1,
$$\align
	\frac1{r} \sum_{s=1}^m \sum_{i \in I(s)} \alpha_i(s) &\leq \frac1{r} 
		\sum_{s=1}^m R_{A(s)}(L) \\
	&= \frac1{r} R_{\sum A(s)}(L) = \frac1{n} \sum_{s=1}^m 
		\sum_{i=1}^n \alpha_i(s).
\endalign$$
Since the first and last terms are equal, we must have $\sum_{i \in I(s)} 
\alpha_i(s) = R_{A(s)}(L)$ for each $s$.  But Corollary 1 to Proposition 1 
implies that, if $\alpha(s)$ consists of $n$ distinct eigenvalues, $L$ is 
unique: there cannot be two or more such $L$ in this intersection. 
\footnote{One can also prove this last part the same way as the first:  If 
($*_\mathscr{I}$) holds with equality, but there are two spaces $L$ and 
$L^\prime$ in $\bigcap_{s=1}^m \Omega_{I(s)}^\circ 
(F_{\sssize{\bullet}}(s))$, the same argument shows that
$$
   	\operatorname{dim}((L+L^\prime)\cap F_i(s)) = 
	\operatorname{dim}(L\cap F_i(s)) + 
	\operatorname{dim}(L^\prime\cap F_i(s)) - 
	\operatorname{dim}((L\cap L^\prime)\cap F_i(s))
$$
for any $s$ and $i$ such that $\alpha_i(s) - \alpha_{i+1}(s) > 0$.  If there is 
an $s$ with all $\alpha_i(s)$ distinct, the display must be an equality for all 
$i$.  This implies, by an easy induction on $i$, that $L\cap F_i(s) = 
L^\prime \cap F_i(s)$ for all $i$, so $L = L^\prime$.}

This argument applies also to the case where the intersection of the 
Schubert classes is not zero or a multiple of the class of a point, for in this 
case the corresponding intersection of Schubert varieties would be infinite.

There is another fundamental fact from geometric invariant theory that 
Klyachko uses to make the link between stability and highest weights.  
Namely, if $\mathscr{L}$ is the restriction to $X$ of the standard line 
bundle $\mathscr{O}(1)$ on 
the projective space $\mathbb{P}(E) = Gr(1,E)$, a point $x$ is semistable 
exactly when some positive power $\mathscr{L}^{\otimes N}$ has a section 
that is invariant by ${SL}(V)$ and does not vanish at $x$.  If each 
$\alpha(s)$ consists of nonnegative integers, the space of sections 
$\Gamma(X,\mathscr{L}^{\otimes N})$ is isomorphic to the dual of the 
tensor product $\otimes_{s=1}^m V(N\alpha(s))$ of representations with 
highest weights $N\alpha(s)$, where, for any $n$-tuple $\alpha  = 
(\alpha_1, \ldots , \alpha_n)$, $N\alpha$  denotes $ (N\alpha_1, \ldots , 
N\alpha_n)$.  The bundle $\mathscr{L}^{\otimes N}$ has a 
nonzero ${SL}(V)$-invariant section if and only if the representation 
$\otimes_sV(N\alpha(s))$ contains the trivial representation of ${SL}(V)$.  
The conditions ($*_\mathscr{I}$) imply that a general point $x$ in $X$ is 
semistable, and therefore $\otimes_s V(N\alpha(s))$ must contain this 
trivial representation.  (For more detail see [Fu3]).  This proves

\proclaim{Proposition 11}  Suppose $\alpha(1), \ldots , \alpha(s)$ are 
partitions of length at most $n$, and \rom{(}$*_\mathscr{I}$\rom{)} holds 
whenever $\prod_s \omega_{I(s)} = \omega_{\{1,\ldots,r\}}$.  
Then there is some 
positive integer $N$ such that $\otimes_s V(N\alpha(s))$ contains the 
trivial representation of ${SL}(V)$. \endproclaim

\head 8.  Coadjoint orbits and moment maps \endhead 

The methods of coadjoint orbits and moment maps from symplectic 
geometry are powerful tools for studying many of the problems 
discussed here.  A full discussion would go beyond the scope of 
this article.  Here we sketch how this approach can be used to 
simplify and unify some of the results obtained in Section 7. An 
excellent account of this method can be found in [Kn].

Let $H$ be the space of Hermitian $n$ by $n$ matrices.  For $\alpha$ 
a weakly decreasing sequence of $n$ real numbers, let $O_\alpha 
\subset H$ be the set of Hermitian matrices with spectrum $\alpha$.  This
is an orbit of the (compact) unitary group $U(n)$ on $H$, acting by 
conjugation.  The Lie algebra of $U(n)$ is the space $\frak{u}(n)$ of 
skew-symmetric matrices.  Identify $H$ with the (real) dual space 
$\frak{u}(n)^*$ by the map which takes $A$ in $H$ to $[B \mapsto 
-\operatorname{tr}(iA\cdot B)]$ in $\frak{u}(n)^*$.  
This identifies $O_\alpha$ with an 
orbit in $\frak{u}(n)^*$, i.e., with a {\bf coadjoint orbit}.

If each $\alpha_i$ is integral, $O_\alpha$ can be 
identified with a partial flag variety $Fl(\alpha)$, which has subspaces of 
dimensions $i$ for each $i$ such that $\alpha_i > \alpha_{i+1}$.  Corresponding 
to $A$ in $O_\alpha$ one has the flag of sums of eigenspaces corresponding 
to $\alpha_1, \ldots , \alpha_i$, for $1 \leq i \leq n$.  As in Section 7, 
there is a natural embedding of $Fl(\alpha)$ in $\mathbb{P}(V(\alpha))$.

Every coadjoint orbit has a natural symplectic structure, as does every 
submanifold of projective space.  It is a basic fact that these symplectic 
structures agree on $O_\alpha \cong Fl(\alpha) \subset \mathbb{P}(V(\alpha))$.

Now suppose $\alpha(1), \ldots , \alpha(m)$ are given, each integral, 
with $\sum_s \sum_i \alpha_i(s) = 0$.  We have closed embeddings 
$$
X = \prod_{s=1}^m O_{\alpha(s)} = \prod_{s=1}^m Fl(\alpha(s)) 
\subset \prod_{s=1}^m \mathbb{P}(V(\alpha(s))) \subset 
\mathbb{P}\left(\otimes_{s=1}^m V(\alpha(s)) \right).
$$
The diagonal action of $U(n)$ on $X$ is a Hamiltonion action, which 
determines a moment mapping from $X$ to $\frak{u}(n)^*$.  Identifying 
$\frak{u}(n)^*$ with $H$, this moment mapping
$$
\mu  \, : \, \, X \to H 
$$
is simply addition: $\mu(A(1), \ldots , A(m)) = A(1) + \ldots + A(m)$.

A basic result of symplectic geometry --- inspired, interestingly, by 
A. Horn's other paper [H1] about eigenvalues of Hermitian matrices --- 
is that the intersection of $\mu(X)$ with a Weyl chamber is a convex 
polytope.  Here it says that the possible eigenvalues $\gamma = 
(\gamma_1 \geq \ldots \geq \gamma_n)$ of a sum $\sum A(s)$ form a 
convex polytope.

In connection with the moment map $\mu$, one has the {\bf symplectic 
reduction} $\mu^{-1}(0)/U(n)$.  This space is in fact homeomorphic to 
the quotient space of $X$ by $G = GL(n,\mathbb{C})$ constructed by 
geometric invariant theory ([MFK],\S8):
$$
\mu^{-1}(0)/U(n) \cong X^{ss}/G \cong \operatorname{Proj}(R), \tag33
$$
where $R = \oplus R_N$, $R_N = \Gamma(X,\mathcal{O}(N))^G = 
\left( \otimes V(N\alpha(s))^* \right)^G$.  From this one sees the 
equivalence of the following three conditions, which say that each 
of the three spaces in (33) is nonempty:
\roster
\item"(i)" there are Hermitian $A(1), \ldots, A(m)$ with eigenvalues 
$\alpha(1), \ldots , \alpha(m)$ such that $\sum_s A(s) = 0$; 
\item"(ii)" a general point (or some point) in $X = \prod Fl(\alpha(s))$ 
is semistable;
\item"(iii)" for infinitely many (or one) $N \geq 1$, $\otimes_s 
V(N\alpha(s))^G \neq 0$.
\endroster

As in Section 7 (cf. [Fu3]), (ii) is equivalent to the condition that 
for general flags $F_{\sssize{\bullet}}(1), \ldots F_{\sssize{\bullet}}(m)$ 
and all subspaces $L$ of $\mathbb{C}^n$,
$$
\sum_{s=1}^m \sum_{i=1}^n \, \alpha_i(s) \, \operatorname{dim} 
\left( L \cap F_i(s) / L \cap F_{i-1}(s) \right) \, \leq \, 0.
$$
If $I(s) = \{ i \in \{1, \ldots , n\} \mid L \cap F_i(s) \neq L \cap 
F_{i-1}(s)\}$, this says that $L$ is in $\bigcap_s \Omega_{I(s)}^\circ 
F_{\sssize{\bullet}}(s)$.  Since the flags are general position, 
(ii) is equivalent 
to 
\roster
\item"(ii${}^\prime$)" $\sum_{s} \sum_{i \in I(s)} \alpha_i(s) \leq 0 \quad 
\text { whenever } \quad \prod_s \omega_{I(s)} \neq 0$.
\endroster

The equivalence of (i), (ii${}^\prime$), and (iii) reproves Propositions 5, 
7, and 11.  For generalizations and variations, see [AW], [BS], [DRW], [He], 
[OS], [SQS], [Ta].

\head 9.  Saturation \endhead

Proposition 11 is close to solving the highest weight problem, except for the 
serious possibility that the integer $N$ that appears might be greater than 
$1$.  In fact, for representations of some classical groups, such integers are 
necessary.  For example, J. Stembridge points out that if $\omega$ is the 
highest weight of the basic spin representation $V(\omega)$ of 
${S@!pin}_9(\mathbb{C})$, then $V(\omega)\otimes V(\omega)$ 
decomposes into four irreducible representations, while 
$V(2\omega)\otimes V(2\omega)$ contains six irreducible representations 
of the form $V(2\eta)$ for some weights $\eta$. Similarly in Schubert 
calculus, P. Pragacz points out that for the Lagrangian Grassmannian of 
3-planes in 6-space, $\sigma_{6@,4@,2}$ occurs in $\sigma_{4@,2}\cdot 
\sigma_{4@,2}$, although $\sigma_{3@,2@,1}$ does not occur in 
$\sigma_{2@,1}\cdot \sigma_{2,@,1}$, cf. [PR].

The key to this is provided by the following result of A. Knutson and T. Tao, 
which solves what is called the Saturation Problem.

\proclaim{Proposition 12 [KT]}  If $\alpha$, $\beta$, and $\gamma$  are a 
triple of partitions, and $c_{N\alpha \,  N\beta}^{\,\,\,N\gamma}  \ne 0$, 
then $c_{\alpha\,\beta}^{\,\,\gamma}  \ne 0$. \endproclaim

Their proof uses a wonderful new combinatorial description of these 
Littlewood-Richardson numbers as the number of geometric figures called 
``honeycombs'' satisfying some conditions.  Their original proof also used 
another new description that they call ``hives.''  Buch [Bu] has given a 
shorter version of their proof, entirely in the language of hives.  For more 
about this Saturation Problem (before it was solved), see A. Zelevinsky's 
article [Z].  Recently H. Derksen and J. Weyman [DW] have given a different  
proof of Proposition 12, using representations of quivers.  It is interesting to 
note that Johnson's unpublished thesis [J] contains a solution of the 
saturation problem when the partitions have lengths at most four.

A. Postnikov and Zelevinsky have pointed out how this result is equivalent to 
one that looks stronger:

\proclaim{Proposition 13}  If $\alpha(1), \ldots , \alpha(m)$ and 
$\beta$  are partitions of length at most $n$, and $N$ is a positive 
integer such that $V(N\beta)$ occurs in the decomposition of 
$\otimes_s V(N\alpha(s))$ as representations of ${GL}_n(\mathbb{C})$, 
then $V(\beta)$ occurs in the decomposition of $\otimes_s 
V(\alpha(s))$. \endproclaim

This follows from that fact that the number of times a representation 
$V(\beta)$ occurs in the decomposition of $\otimes_s V(\alpha(s))$ is 
equal to a Littlewood-Richardson coefficient 
$c_{\alpha\,\beta}^{\,\,\gamma}$, where $\alpha$  and $\gamma$  are 
obtained as follows.  Arrange the Young diagrams of $\alpha(1), \alpha(2), 
\ldots , \alpha(m)$ corner to corner from lower left to upper right, so that 
the upper right corner of the diagram of $\alpha(i)$ just touches the lower 
left corner of $\alpha(i+1)$; these form a skew diagram 
$\gamma \smallsetminus \alpha$.  The reason for the equality is that both 
numbers are the numbers of skew tableaux on the shape $\gamma 
\smallsetminus \alpha$  whose rectification is a given tableau of shape 
$\beta$  (cf. [Fu2], \S5.1).  Proposition 13 can also be proved directly, using 
the methods of [KT] or [DW].

Note that the corresponding result follows for representations of 
${SL}_n(\mathbb{C})$, since representations of ${SL}_n(\mathbb{C})$ 
differ from those of ${GL}_n(\mathbb{C})$ only by powers of the 
determinant representation $\bigwedge^n\mathbb{C}^n$ (cf. [FH], \S15).  

\proclaim{Proposition 14}  Let $\alpha(1), \ldots , \alpha(m)$ be partitions 
of lengths at most $n$, and suppose $c = \frac1{n} \sum_{s=1}^m \sum 
\alpha_i(s)$ is an integer.  Then the representation $\otimes_{s=1}^m 
V(\alpha(s))$ of ${GL}_n(\mathbb{C})$ contains the representation 
$(\bigwedge^n\mathbb{C}^n)^{\otimes c}$ if and only if 
\rom{(}$*_\mathscr{I}$\rom{)} holds for all $\mathscr{I} \in S_r^n(m)$ 
and all $r < n$. \endproclaim

These conditions are equivalent to the condition that the Schur polynomial 
$s_{(c^n)}$ occurs in the product $\prod_{s=1}^m s_{\alpha(s)}$.

As to geometric invariant theory itself, the ring $\oplus_{N\geq 0}  
\Gamma(X,\mathscr{L}^{\otimes N})^{{SL}(V)}$ is the homogeneous 
coordinate ring of the GIT quotient variety $X/\!\!/{SL}(V)$ 
of $X$ by ${SL}(V)$; that is, $X/\!\!/{SL}(V) = \operatorname{Proj}(R)$.
One may ask whether this graded ring $R$ is generated by its 
homogeneous part of degree $1$.  Unfortunately, this is not always true, 
even for these special varieties that are products of partial flag varieties:  

\example{Example 5}  For $n = 3$, let $X$ be the Cartesian product of $6$ 
copies of $\mathbb{P}^2 = \mathbb{P}(V)$, and $\mathscr{L} = 
\mathscr{O}(1,1,1,1,1,1)$ the tensor product of the pullbacks of the line 
bundles $\mathscr{O}(1)$ on the factors.  The homogeneous part 
$\Gamma(X,\mathscr{O}(2,2,2,2,2,2))^{{SL}_3(\mathbb{C})}$ of degree 
$2$ has dimension $16$, while that of degree $1$ has dimension $5$, and 
$16 > \binom{5+1}{2}$.  Equivalently, if $\alpha  = (5,4,3,2,1)$, $\beta  = 
(2,2,2)$, and $\gamma  = (6,5,4,3,2,1)$, then 
$c_{\alpha\,\beta}^{\,\,\gamma}  = 5$, while $c_{2\alpha \, 
2\beta}^{\,\,\,2\gamma}  = 16 > \binom{c_{\alpha\,\beta}^{\,\,\gamma} 
+1}{2}$.  Buch has verified by computer that the smallest $N \cdot \sum 
\gamma_i$ for which there is an $\alpha$, $\beta$, $\gamma$  and $N$ 
for which $c_{N\alpha \, N\beta}^{\,\,\,N\gamma}  > 
\binom{c_{\alpha\,\beta}^{\,\,\gamma} + N-1}{N}$ is for $N \cdot \sum 
\gamma_i = 30$, with $\alpha  = (4,3,2,1)$, $\beta  = (2,2,1)$, $\gamma  
= (5,4,3,2,1)$, and $N = 2$, where again $c_{\alpha\,\beta}^{\,\,\gamma}  
= 5$, and $c_{2\alpha \,  2\beta}^{\,\,\,2\gamma}  = 16$. 
For a geometric description of $(\mathbb{P}^2)^6/\!\!/SL(3,\mathbb{C})$, 
see [DO], \S{VII.3}.  \endexample

A. Okounkov speculates that the function 
$f(N) = c_{N\alpha \, N\beta}^{\,\,\,N\gamma}$ may be log concave, i.e., that 
$$
f(N)^2 \geq f(N-1)\cdot f(N+1)
$$
for all $N \geq 1$.

For a general product $X$ of partial flag varieties, is there a reasonable 
criterion to tell when this algebra is generated by its part in degree $1$, or 
to find an $N$ such that the algebra 
$\oplus_{k\geq 0}  \Gamma(X,\mathscr{L}^{\otimes kN})^{{SL}(V)}$ is 
generated by $\Gamma(X,\mathscr{L}^{\otimes N})^{{SL}(V)}$?

Knutson, Tao, and Woodward [KTW] 
 have also announced a proof of a conjecture we had made 
that $c_{N\alpha \,  N\beta}^{\,\,\,N\gamma}  = 1$ if and only if 
$c_{\alpha\,\beta}^{\,\,\gamma}  = 1$.  Buch [Bu] deduced from this that, 
for $n\geq  3$, if each of the partitions in $(\alpha,\beta,\gamma)$ consists 
of $n$ distinct integers, and if each of the inequalities ($*_{IJK}$) holds 
strictly, then $c_{\alpha\,\beta}^{\,\,\gamma}$  must be at least $2$.

\head 10.  Proofs of the theorems \endhead

Most of the theorems in Sections 1--4 follow easily from the propositions 
proved in Sections 6, 7, and 9.  Let us 
start by summarizing the main results, but 
generalized to an arbitrary number of factors.  Alternative versions are stated 
in brackets.
\subhead 10.1. The main theorem \endsubhead

\proclaim{Theorem 17}  Let $\alpha(1), \ldots , \alpha(m)$ and $\gamma$  
be weakly decreasing sequences of $n$ real numbers, with $\sum_{i=1}^n 
\gamma_i = \sum_{s=1}^m \sum_{i=1}^n \alpha_i(s)$.  The following are 
equivalent:
\roster
\item"(1)"  $\sum_{k \in K} \gamma_k \leq \sum_{s=1}^m \sum_{i \in 
I(s)} \alpha_i(s)$ is satisfied 
for all subsets $I(1), \ldots , I(m)$, $K$ of the same 
cardinality $r$ of $\{1, \ldots , n\}$ such that 
$\sigma_{\lambda(K)}$ occurs in $\prod_{s=1}^m \sigma_{\lambda(I(s))}$ 
 \rom{[}with coefficient $1$\rom{]} 
in $H^*(Gr(r,\mathbb{C}^n))$, and all $r < n$.
\item"(2)"  $\alpha(1), \ldots , \alpha(m)$ and $\gamma$  are eigenvalues 
of Hermitian 
matrices $A(1), \ldots , A(m)$ 
and $C$, such that $C = \sum_{s=1}^m A(s)$.
\endroster
If each of $\alpha(1), \ldots , \alpha(m)$ and $\gamma$  consists of 
integers, (1) and (2) are equivalent to 
\roster
\item"(3)"  The representation $V(\gamma)$ of ${GL}_n(\mathbb{C})$ 
occurs in $\otimes_{s=1}^m V(\alpha(s))$.
\endroster
If each of  $\alpha(1), \ldots , \alpha(m)$ and $\gamma$  is a partition, 
these are equivalent to 
\roster
\item"(4)"  For some 
\rom{[}every\rom{]} discrete valuation ring $R$, there exist $n$ by 
$n$ matrices $A(1), \ldots , A(m)$ and $C$ with invariant factors 
$\alpha(1), \ldots , \alpha(m)$ and $\gamma$  such that $C = A(1)\cdot  
\ldots \cdot A(m)$. 
\endroster
If, in addition, $q$ is an integer at least as large as each  
$\alpha_1(s)$ and $\gamma_1$, these are equivalent to
\roster
\item"(5)" The class $\sigma_\gamma$  occurs with positive coefficient in 
the product 
$\prod_{s=1}^m \sigma_{\alpha(s)}$ in $H^*(Gr(n,\mathbb{C}^{n+q}))$.
\endroster
\endproclaim
We show how to deduce this theorem from the propositions in the 
preceding sections.  As in the proof of Proposition 2, $\sigma_{\lambda 
(K)}$ occurs in $\prod_{s=1}^m \sigma_{\lambda(I(s))}$ [with coefficient 
$1$] if and only if $\prod_{s=1}^m\omega_{I(s)^\prime} \cdot \omega_K$  
is a nonzero multiple [equal to $1$] of $\omega_{\{1,\ldots,r\}}$.  
Let $\beta(s) = (-
\alpha_n(s), \ldots , -\alpha_1(s))$, for $1 \leq s \leq m$, and let 
$\beta(m+1) = \gamma$. 
The condition (2) is equivalent to the existence of Hermitian 
matrices $B(s)$ with eigenvalues $\beta(s)$, $1 \leq s \leq m+1$, whose 
sum is zero (as one sees by setting $A(s) = -B(s)$ and $C = B(m+1)$).  From 
Proposition 7 [and Proposition 9], (2) is equivalent to the inequalities 
$\sum_{s=1}^m \sum_{i \in I(s)^\prime} \beta_i(s) + \sum_{k \in 
\mathscr{K}} \gamma_k \leq 0$ whenever 
$\prod_{s=1}^m\omega_{I(s)^\prime} \cdot \omega_K$ is a nonzero 
multiple [equal to $1$] of $\omega_{\{1,\ldots,r\}}$.  
This proves the equivalence of 
(1) and (2).

Condition (3) is equivalent to $\otimes_{s=1}^m V(\alpha(i))^*\otimes 
V(\gamma)$ containing the trivial representation of ${GL}_n(\mathbb{C})$, 
where $V(\alpha(s))^*$ is the dual representation to $V(\alpha(s))$, which 
has highest weight $(-\alpha_n(s), \ldots , -\alpha_1(s))$.  So the 
equivalence of (3) with (1) follows from Propositions 11 and 13.  

Condition (3) is equivalent to the assertion that the Schur polynomial 
$s_\gamma$  occurs in the product $\prod_{s=1}^m s_{\alpha(s)}$.  The 
equivalence of this with (4) is seen by induction on $m$, the case $m = 2$ 
having been seen in Section 6.  Indeed, $s_\gamma$  occurs in  
$\prod_{s=1}^m s_{\alpha(s)}$ if and only if there is 
some $\beta$  such that $s_\beta$  occurs in $\prod_{s=1}^{m-1} 
s_{\alpha(s)}$ and $s_\gamma$  occurs in $s_\beta \cdot s_{\alpha(m)}$.  
By induction, the first corresponds to a product $A(1)\cdot  \ldots \cdot 
A(m-1) = B$ with invariant factors $\alpha(1), \ldots , \alpha(m-1)$ and 
$\beta$, and the second to a product $B\cdot A(m) = C$ with invariant 
factors $\beta$, $\alpha(m)$, and $\gamma$.

That (3) and (5) are equivalent follows as before from the fact that 
there is a homomorphism from the ring of symmetric polynomials onto the 
cohomology of the Grassmann manifold taking $s_\alpha$  to 
$\sigma_\alpha$  for all $\alpha$.

\subhead 10.2. Horn's conjecture \endsubhead
We turn to the proof of Theorem 12, which implies Horn's conjecture.  We 
start with a general combinatorial lemma. 

\proclaim{Lemma 4} (i) Let $I = \{i_1 < \ldots < i_r\} \subset \{1, 
\ldots ,n\}$, and let $\lambda  = (i_r - r, \ldots , i_1 - 1)$ be the 
corresponding partition.  Let $I' = \{i \in \{1, \ldots , n\} \,|\, 
n+1-i \notin I\}$.  Then the partition corresponding to $I'$ is 
the conjugate partition $\lambda'$  to $\lambda$. \par 
(ii) If $\lambda$, $\mu$, and $\nu$ are partitions, with conjugates 
$\lambda'$, $\mu'$, and $\nu'$, then  
$c_{\lambda \, \mu}^{\,\,\nu} = c_{{\lambda'} \,  
{\mu'}}^{\,\,{\nu'}}$.
\endproclaim
The first is an elementary combinatorial fact 
(cf. [Mac] \S I(1.7)) that can be seen visually 
by identifying $\lambda$  with its Young diagram, which is traced out by a 
sequence of vertical steps going from lower left to upper right of the $r$ by 
$n-r$ rectangle of boxes.  If $\lambda = \lambda(I)$, then $I$ specifies 
which of these $n$ steps are vertical, and its complement specifies which 
are horizontal, so $I'$ specifies the vertical steps for 
$\lambda'$.  The second assertion is 
not obvious from the definition of Littlewood-Richardson coefficients, but 
follows from other descriptions, cf. [Fu2] \S5.1.

Now we prove that $T_r^n = S_r^n$ for all $r < n$, by induction on $r$.  It 
follows from the definition that $T_1^n = S_1^n = U_1^n$.  Let $(I,J,K)$ 
be a triple in $U_r^n$, with $r > 1$, and let $(\lambda ,\mu,\nu)$ be the 
corresponding triple of partitions, each of length at most $r$.  The 
condition for $(I,J,K)$ to be in $T_r^n$ is that, for all $p < r$ and all 
$(F,G,H) \in T_p^r$, 
$$
	\sum_{f \in F} i_f +\sum_{g \in G} j_g  \leq \sum_{h \in H} k_h + 
		p(p+1)/2. 
$$
Since $i_a - a = \lambda_{r+1-a}$, and similarly for $J$ and $K$, this is 
equivalent to the inequality
$$
	\sum_{f \in F} \lambda_{r+1-f} +\sum_{g \in G} \mu_{r+1-g}  \leq 
		\sum_{h \in H} \nu_{r+1-h}.   
$$
The assumption that $(I,J,K)$ is in $U_r^n$ says that the triple 
$(\lambda,\mu,\nu)$ satisfies ($*$), i.e., $\sum \lambda_i + \sum \mu_i = 
\sum \nu_i$.  The preceding inequality is therefore equivalent to 
$$
	\sum_{f \in F'} \lambda_f +\sum_{g\in G'} 
		\mu_g \geq  \sum_{h\in H'} \nu_h ,
$$
where $F' = \{f \in \{1, \ldots , r\} \,|\, r+1-f \notin F\}$, and 
$G'$ and $H'$ defined similarly.
By induction, we know that $T_p^r = S_p^r$, and by Lemma 4 we know 
that $(F,G,H)$ is in $S_p^r$ if and only if 
$(F',G',H')$ is in $S_{r-p}^r$.  Hence the 
inequalities in question are equivalent to the inequalities 
$$
\sum_{h\in H} \nu_h \leq \sum_{f \in F} \lambda_f +\sum_{g \in G} 
\mu_g 
$$
for all $q = r-p < r$ and all $(F,G,H)$ in $S_q^r$.  But by the equivalence of 
(3) and (5) of Theorem 17 (for the case $m = 2$), we know that 
the validity of these 
inequalities ($*_{FGH}$) for all $(F,G,H)$ in $S_p^r$ for all $p < r$ is  
equivalent to the assertion that $c_{\lambda \, \mu}^{\,\,\nu}$ is not zero. 
This means that $(I,J,K)$ is in $S_r^n$.

The same idea proves the following version for an arbitrary number $m$ of 
factors, but in a symmetric form where all factors are treated equally (so 
without a duality argument).  Let $U_r^n(m)$ be the set of $m$-tuples 
$I(1), \ldots , I(m)$ of subsets of cardinality $r$ in $\{1, \ldots ,n\}$ such 
that 
$$
\sum_{s=1}^m \sum_{i\in I(s)} i = (m-1)r(n-r) + mr(r+1)/2. \tag34
$$
Define $S_r^n(m) \subset U_r^n(m)$ to be those $m$-tuples such that the product 
$\prod \omega_{I(s)}$ of the corresponding Schubert varieties is a nonzero 
multiple of the class of a point in $Gr(r,\mathbb{C}^n)$, as 
in Section 6.  Define $T_r^n(m) \subset U_r^n(m)$ inductively, to be the 
set such that, for all $p < r$ and all $\mathscr{F} = (F(1), \ldots , F(m))$ in 
$T_p^r(m)$,
$$
\sum_{s=1}^m\sum_{f\in F(s)} i_f(s) \geq  (m-1)p(n-r) + mp(p+1)/2. 
\tag35
$$
\proclaim{Theorem 18}  $S_r^n(m) = T_r^n(m)$ for all positive integers 
$m$, $r$, and $n$ with $r < n$. \endproclaim
This is easier to prove.  For $r = 1$, $T_1^n(m) = S_1^n(m) = U_1^n(m)$.  
For $r > 1$, let $\mathscr{I} = (I(1), \ldots ,I(m))$ be in $U_r^n(m)$, 
and define $\lambda(s)$ so that $\sigma_{\lambda(s)} = \omega_{I(s)}$, i.e., 
$i_j(s) = n-r+j-\lambda_j(s)$ for $1 \leq j \leq r$.  Then $\mathscr{I}$ is 
in $T_r^n(m)$ if and only if $\sum_s \sum_{f \in F(s)}(i_f(s) - f) \geq (m-
1)p(n-r)$ for all $\mathscr{F}$ in $T_p^r(m)$ and all $p < r$.  This is 
equivalent to the inequality $\sum_s \sum_{f \in F(s)} \lambda_f(s) \leq 
p(n-r)$, or $\frac1{p} \sum_s \sum_{f \in F(s)} \lambda_f(s) \leq n-r$ for 
all $\mathscr{F}$ in $T_p^r
(m)$ and all $p < r$.  By Proposition 14, and the inductive fact that 
$T_p^r(m) = S_p^r(m)$, this is equivalent to saying that $\prod_s 
\sigma_{\lambda(s)} = d \, \sigma_\rho$ for some $ d \ne 0$, where $\rho 
= (n-r)^r$.  This says that $\mathscr{I}$ is in $S_r^n(m)$. 

By proving Theorems 17 and 18, we have completed the proofs of Theorems 1
and 2 in their general 
forms for an arbitrary number of factors.  Theorem 5, also for any number of 
factors, follows similarly from Proposition 8. 

\subhead 10.3.  Compact operators \endsubhead  To prove Theorem 6, for any 
number of compact self-adjoint operators $A(1), \ldots , A(m)$ and $C = 
\sum_{s=1}^m A(s)$ on a Hilbert space $H$, we want to show that all the 
inequalities of (1) of Theorem 17 are valid.  Choose a finite-dimensional 
subspace $V$ of $H$ that contains eigenvectors corresponding to the $n$ 
largest eigenvalues of each of the operators under consideration.  It then 
suffices to apply the finite-dimensional result to the operators $A(1)_V, 
\ldots , A(m)_V$, and $C_V$;  here, for any operator $A$, $ A_V$  
denotes the composite $V \to H \to H \to V$, where the first map is the 
inclusion, the second $A$, and the third is orthogonal projection.  The same 
argument as in Corollary 2 to Proposition 1 shows that the $k^{\text{th}}$ 
largest eigenvalue of $A_V$ is no larger than the $k^{\text{th}}$ largest 
eigenvalue of $A$.  The operators $A(s)_V$ and $C_V$ will therefore 
have the same 
$n$ largest eigenvalues as $A(s)$ and $C$.  For more on eigenvalues of 
compact operators, see [W], [Zw], [Fr], and [R].

\subhead 10.4.  Invariant factors \endsubhead
We have seen the general version of Theorem 7 in Theorem 17.  
Theorem 8 follows from 
Theorem 7 by localization.  This is easiest to see in the version about short 
exact sequences of modules of finite length.  To specify such an exact 
sequence over a principal ideal domain $R$ is equivalent to specifying a 
short exact sequence of modules of finite length over a finite number of 
localizations $R_{\mathfrak{p}}$ of $R$, each of which is a discrete 
valuation ring.  

Theorem 9 follows easily from Theorem 8.  Indeed, given $C$ (i.e., $X$), 
there is an exact sequence 
$$
	0 \to \mathcal{A} \to \mathcal{C} \to \mathcal{B} \to 0, 
$$
where $\mathcal{A}$, $\mathcal{B}$, and $\mathcal{C}$ are the cokernels 
of $A$, $B$, and $C$, so the necessity follows from Theorem 8.  Conversely, 
given $A$ and $B$ with invariants $a_1 \subset \ldots \subset a_p$ and  
$b_1 \subset \ldots \subset b_q$, and possible invariants $c_1 \subset 
\ldots \subset c_n$ satisfying the conditions displayed in Theorem 9, 
Theorem 8 implies that there is such an exact sequence, with 
$\mathcal{A}$ and $\mathcal{B}$ the cokernels of $A$ and $B$.  The 
standard argument used to resolve a short exact sequence of modules by a 
long exact sequence of free modules says that the dotted arrows in the 
following diagram can be filled in: 
$$\xymatrix{
R^p \ar[r]^i \ar[d]_A & R^p\oplus R^q  \ar[r]^j \ar@{.>}[d]  & R^q \ar[d]^B 
\\
R^p \ar[r]^i \ar[d] & R^p\oplus R^q \ar[r]^j \ar@{.>}[d] & R^q \ar[d] \\
\mathcal{A}  \ar[r] &  \mathcal{C}     \ar[r]  & \mathcal{B}
}
$$
Here $i$ and $j$ are the canonical embeddings and projections: $i(v) = 
(v,0)$, $j(v,w) = w$.  This produces a matrix $C$ of the required form with 
cokernel $\mathcal{C}$.

To deduce Theorem 10 from Theorem 9, we use the standard fact that an 
endomorphism of a finite dimensional vector space over a field $F$, up to 
similarity, is the same as a finitely generated torsion module over $R = 
F[T]$, up to isomorphism.  The invariant factors are the same for the 
endomorphism and the $R$-module, and they determine the 
endomorphism up to similarity. 

\subhead 10.5.  Highest weights \endsubhead
Theorems 11 and 14 are special cases of Theorem 17 together with 
Theorem 12. Theorem 13 follows from Proposition 9, together with the fact 
that $(I,J,K)$ is in $S_r^n$ (resp. $R_r^n$) if and only if 
$(I^\prime,J^\prime,K)$ is in $S_r^n(3)$ (resp. $R_r^n(3)$).
Here $R_r^n(m)$ is defined to be the subset of $S_r^n(m)$ such that 
the product $\prod \omega_{I(s)}$ is equal to the class of a point.

\subhead 10.6.  Singular values \endsubhead
We next prove the theorems about singular values.  For the additive version 
in Theorem 15, the fact that all the inequalities ($**_{IJK}$) must be 
satisfied follows immediately from Theorem 1, since if $C = A+B$ is a sum of 
arbitrary matrices, then  
$$
	\pmatrix 0 & C \\ C^* & 0 \endpmatrix
	\,\, = \,\, \pmatrix 0 & A \\ A^* & 0 \endpmatrix \, + \, 
		\pmatrix 0 & B \\ B^* & 0 \endpmatrix     
$$
is a sum of Hermitian matrices.  (In fact, if one uses the set $S_r^n$ instead 
of the set $T_r^n$, the necessity follows from the simple Proposition 2.)  
The full theorem was recently proved by O'Shea and Sjamaar [OS], \S9, (also 
in terms of the sets $S_r^n$), using moment maps, although they do not 
mention singular values.  Unlike the other theorems described here, we do 
not know another proof for this converse.  The same arguments prove the 
generalization of Theorem 15 to any number of factors.

The easy half of the multiplicative version (Theorem 16) can be proved 
by using a variation of the 
Hersch-Zwahlen lemma, which is based on an idea of Amir-Mo\'ez [AM], cf. [TT1].  
First we require a multiplicative analogue of the Rayleigh trace.  Suppose  
$A$  is a positive semidefinite Hermitian $n$ by $n$ matrix.  Let $\alpha_1 
\geq  \ldots  \geq  \alpha_n \geq  0$ be the eigenvalues of $A$, and let 
$v_1, \ldots , v_n$ be corresponding orthonormal eigenvectors.  Let 
$F_{\sssize{\bullet}}(A)$ be the flag with $F_k(A)$ spanned by $v_1, \ldots 
, v_k$.  If $L \subset \mathbb{C}^n$ is an $r$-dimensional subspace, define 
$D_A(L)$ to be the determinant of the matrix $(A\,u_i,u_j)_{1\leq i,j\leq 
r}$, where $u_1, \ldots , u_r$ is an orthonormal basis of $L$.  Equivalently, 
$D_A(L)$ is the determinant of the mapping $A_L$ of Corollary 3 to 
Proposition 1.   

\proclaim{Proposition 15}  For any subset $I \subset \{1, \ldots , n\}$ of 
cardinality $r$,
$$
	\prod_{i \in I} \alpha_i = 
	\operatornamewithlimits{Min}_{L \in 
\Omega_I(F_{\sssize{\bullet}}(A))} D_A(L).
$$ \endproclaim
To prove this, write $u_i = \sum_{k=1}^n x_{ik}v_k$.   Then $(A\,u_i,u_j)$ 
is the $i,j$ entry of the matrix $X \, D \, X^*$, where $X$ is the $r$ by $n$ 
matrix $(x_{ij})$, and $D$ is the diagonal matrix with diagonal entries 
$\alpha_1, \ldots , \alpha_n$.  Therefore 
$$
   	D_A(L)  =  \sum_J \det(X^J) \cdot \prod_{j \in J} \alpha_j  \cdot   
		\det(X^*_J)  =  \sum_J |\det(X^J)|^2  \cdot   	
		\prod_{j \in J} \alpha_j,
$$
where the sum is over all subsets $J$ of $\{1, \ldots , n\}$ of cardinality 
$r$, $X^J$ is the submatrix of $X$ with columns labeled by $J$, and 
$X^*_J$ is the submatrix of $X^*$ with rows labeled by $J$.  Take $u_p$ in 
$L\cap F_{i_p}$ successively as in the proof of Proposition 1.  The 
$p^{\text{th}}$ row of the matrix $X$ ends in column $i_p$, so $\det(X^J) 
= 0$ unless $j_p \leq i_p$ for all $p$.  For such $J$,  $\prod_{j \in J} 
\alpha_j  \geq  \prod_{i \in I} \alpha_i$.  Moreover, $\sum_J 
|\operatorname{det}(X^J)|^2 \geq  \det(XX^*) = 1$.  These inequalities 
imply that $D_A(L) \geq  \prod_{i \in I} \alpha_i$.  As before, equality is 
achieved by taking $L$ to be spanned by those vectors $v_i$ for $i$ in $I$.

\proclaim{Proposition 16}  Let $C = A(1) \cdot \ldots \cdot A(m)$ be a 
product of $m$ arbitrary $n$ by $n$ matrices.  Suppose $L_1, \ldots , 
L_m$ are $r$-dimensional subspaces of $\mathbb{C}^n$ with $A(s)(L_s) 
\subset L_{s-1}$ for $2 \leq s \leq m$.  Then 
$$
	D_{C^*C}(L_m) = \prod_{s=1}^m D_{A(s)^*A(s)} (L_s).
$$ \endproclaim
We first do the case where $m = 2$.  Suppose $C = A\cdot B$, and $L$ and 
$M$ are $r$-dimensional linear subspaces with $B(M) \subset L$.  We must 
show that $D_{C^*C}(M) = D_{A^*A}(L) \cdot D_{B^*B}(M)$.  Let $u_1, 
\ldots , u_r$ be an orthonormal basis for $L$, and $v_1, \ldots , v_r$ an 
orthonormal basis for $M$.  Write $B(v_i) = \sum x_{ij}u_j$, and let $X = 
(x_{ij})$.  Then
$$\align
	D_{C^*C}(M) &= \det((B^*A^*ABv_i,v_j)) = \det((A^*ABv_i,Bv_j)) \\
		&= \det(X) \cdot \det((A^*Au_i,u_j)) \cdot \det(X^*) 
		= \det(XX^*) \cdot \det((Au_i,Au_j)) \\
		&=  \det(XX^*) \cdot D_{A^*A}(L).
\endalign$$
But $D_{B^*B}(M) = \det((Bv_i,Bv_j)) = \det(\sum x_{ik} \overline{x_{jl}} 
(u_k,u_l)) = \det(XX^*)$, which proves the claim.  The general case is 
proved from this by induction, setting $A = A(1)\cdot  \ldots \cdot A(m-
1)$, $B = A(m)$, $M = L_m$, and $L = L_{m-1}$.

\proclaim{Theorem 19}  Let $A(1), \ldots , A(m)$ be $n$ by $n$ matrices 
such that $A(1) \cdot  \ldots \cdot A(m) = 1$.  Let $a_1(s) \geq  \ldots  
\geq  a_n(s) > 0$ be the singular values of $A(s)$, $1 \leq s \leq m$.  
Suppose $I(1), \ldots , I(m)$ are subsets of cardinality $r$ of $\{1, \ldots 
,n\}$ such that the product $\prod \sigma_{I(s)}$ of the corresponding 
Schubert classes is not zero.  Then 
$$
	\prod_{s=1}^m \prod_{i \in I(s)} a_i(s) \leq 1.
$$\endproclaim
Note that $a_1(s)^2 \geq  \ldots  \geq  a_n(s)^2$ are the eigenvalues of 
$A(s)^*A(s)$.  To prove the theorem, take flags $F_{\sssize{\bullet}}(s)$ 
corresponding to orthonormal bases of eigenvectors for $A(s)$.  By 
Proposition 15, $D_{A(s)^*A(s)}(L_s) \geq  \prod_{i \in I(s)} a_i(s)^2$ if 
$L_s$ is in $\Omega_{I(s)}(F_{\sssize{\bullet}}(s))$.  By the hypothesis on 
the nonvanishing of the product of Schubert classes, we may find a subspace 
$L$ in the intersection 
$$
	\Omega_{I(m)}(F_{\sssize{\bullet}}(m)) 
\,\, {\tsize \bigcap}  \bigcap_{s=1}^{m-1} 
	\Omega_{I(s)}(A(s+1) \cdot \ldots \cdot A(m))^{-
1}(F_{\sssize{\bullet}}(s)). 
$$
Set $L_s = A(s+1) \cdot  \ldots \cdot A(m)(L)$ for $1 \leq s < m$.  Then 
$A(s)(L_s) \subset L_{s-1}$, and $L_s$ is in 
$\Omega_{I(s)}(F_{\sssize{\bullet}}(s))$.  By the two preceding 
propositions, 
$$
	1 = D_1(L) = \prod_{s=1}^m D_{A(s)^*A(s)}(L_s) \geq  
		\prod_{s=1}^m \prod_{i \in I(s)} a_i(s)^2.
$$
\proclaim{Corollary}  Let $A(1), \ldots , A(m)$ be $n$ by $n$ matrices, and 
let $ C = A(1) \cdot \ldots \cdot A(m)$.  Let  $a(s)$  be the singular values 
of $A(s)$, and $c$ the singular values of $C$.  Then for any $I(1), \ldots , 
I(m)$ and $K$ subsets of  $\{1, \ldots ,n\}$ of cardinality $r$, such that 
$\sigma_{\lambda(K)}$ appears in $\prod_{s=1}^m 
\sigma_{\lambda(I(s))}$,
$$
	\prod_{k \in K} c_k \leq \prod_{s=1}^m \prod_{i \in I(s)} a_i(s).
$$\endproclaim
The proof is similar to that used in the additive case for eigenvalues:  By a 
continuity argument, as in the proof of Proposition 7, it suffices to consider 
the case where all the matrices are invertible.  If $a_1 \geq  \ldots  \geq  
a_n > 0$ are the singular values of $A$, then $1/a_n \geq  \ldots  \geq  
1/a_1$ are the singular values of $A^{-1}$.  Apply the theorem to $C \cdot 
A(m)^{-1} \cdot \ldots \cdot A(1)^{-1} = 1$, using the sets $K$ and 
$I(m)^\prime, \ldots , I(1)^\prime$, where $I^\prime$ denotes $\{n+1-i 
\,|\, i \in I\}$, exactly as in the proof of Theorem 17.  

In another remarkable paper [Kl2], Klyachko has recently 
proved the converse, and the corollary 
itself, which specializes to Theorem 16 for $m = 2$.  He deduces this 
from a computation of probability densities for random walks on the Lie 
algebra $\frak{su}(n)$ and the homogeneous space $SL(n,\mathbb{C})/SU(n)$. 
Such densities for Euclidean 3-space were computed by Lord Rayleigh.

\subhead 10.7.  Real symmetric and Hermitian quaternionic matrices 
\endsubhead We prove the general versions of Theorems 3 and 4, Proposition 
7, and Theorem 17 for real symmetric and Hermitian (self-conjugate) 
quaternionic matrices.

\proclaim{Theorem 20} The eigenvalues $\alpha(1), \ldots, \alpha(m)$, and 
$\gamma$ that occur for complex Hermitian matrices $A(1), \ldots, A(m)$ and 
$C = A(1) + \ldots + A(m)$ are the same as the eigenvalues that occur for 
real symmetric or for Hermitian quaternionic matrices. \endproclaim

For real symmetric matrices, it suffices to show that Proposition 7 is 
true for real symmetric matrices.  The argument is similar to the complex 
case, using the fact that if $V$ is the complexification of a real vector 
space, then flags that are complexifications of real flags are Zariski 
dense in all flags.  From this it follows that the Hermitian matrices 
constructed in Proposition 7 can all be taken to be real symmetric 
matrices. (See [Fu3] for details.)

A Hermitian quaternionic matrix is an $n$ by $n$ matrix A with entries in 
the quaternions $\mathbb{H}$ which is equal to the transpose of its 
conjugate.  Such a matrix acts on the left on the (right) 
$\mathbb{H}$-space
$\mathbb{H}^n$ of column vectors. Such a matrix again has all its 
eigenvalues real, and $\mathbb{H}^n$ has an orthonormal basis of 
eigenvectors.  In this case, one needs only prove that all the conditions 
($*_\mathscr{I}$) are necessary.  This follows exactly as in the proofs of 
Propositions 2 and 7, together with the observation that Schubert calculus 
for the quaternionic Grassmannian $Gr(r,\mathbb{H}^n)$ of (right) 
subspoaces of $\mathbb{H}^n$ is exactly the same (but with all degrees 
doubled) as that for complex Grassmannians $Gr(r,\mathbb{C}^n)$. This is 
proved in [PR].

Steinberg asks if one can characterize the eigenvalues of real 
skew-symmetric matrices $A$, $B$, and $C = A+B$.  For any real matrix 
$A$, the matrix $\left(\smallmatrix 0 & A \\ -A^t & 0
\endsmallmatrix\right)$ is skew-symmetric, so this may be related to 
the problem of characterizing singular values of real matrices 
$A$, $B$, and $C = A+B$.
 
\head 11.  Final remarks \endhead

It would be interesting to find more direct relations between the subjects 
addressed here, that would give better explanations of why the 
questions in each subject have the 
same answers.  It follows from the main theorems in this 
article that in each of these 
subjects, the question of which triples $(\alpha,\beta,\gamma)$ occur for a 
given $n$ is completely determined by knowing the answer to the same 
question for all $r < n$.  
It is a challenge to find direct explanations for this 
in any of these situations.  For example, is there a representation-theoretic 
explanation of the fact that knowing which irreducible representations of 
${GL}_n(\mathbb{C})$ occur in tensor products is determined by the 
answer to this same question for ${GL}_r(\mathbb{C})$ for $r < n$?  A similar 
question can be asked about the Schubert calculus for 
$Gr(n,\mathbb{C}^{n+q})$ being determined by the Schubert calculus for all 
$Gr(r,\mathbb{C}^n)$ for $r < n$ (at least the part about which Schubert 
classes appear in a product). 
One has a similar mystery for the invariant factors problems. 
For the eigenvalue problem, Theorem 5 and Proposition 8 at least produce 
some smaller matrices.

Agnihotri and Woodward [AW] and Belkale [Be] have proved an analogue of 
the Klyachko theorem for unitary matrices.  They characterize the possible 
eigenvalues of unitary matrices $A(1), \ldots , A(m)$ whose product is the 
identity.  Instead of being controlled by Schubert calculus as in the 
Hermitian eigenvalue problem, it is controlled by {\em quantum Schubert 
calculus}. Although there are algorithms for computing in this quantum 
Schubert calculus, there are not yet true analogues of the 
Littlewood-Richardson
rule, nor are there useful criteria for the nonvanishing of such 
products. 

One would like to see good analogues of these theorems for other Lie groups.  
We know that at least the saturation problem must be modified, but it seems 
reasonable to hope for sharp analogues of many of the theorems, at least for 
the other classical groups.  The best results along these lines so far have 
been achieved by use of moment maps and coadjoint orbits, cf. [AW], [BS], 
[He], [OS].  Knutson and Tau end their paper [KT] with a precise conjecture. 

One of the relations that comes from this work is that between 
representations of the symmetric group and Hermitian eigenvalues.  There 
is a deep relation, conjectured by Baik, Deift, and Johansson [BDJ], and 
proved by Okounkov [O].  Any permutation $w$ in $S_n$, if written as a 
sequence of its values $w_1, \ldots , w_n$, determines a partition 
$\lambda(w)$ of $n$: $\lambda_1(w)$ is the length of the longest 
increasing subsequence of $w$, $\lambda_1(w) + \lambda_2(w)$ is the maximum 
sum 
of the lengths of the two disjoint increasing subsequences, and so 
on.  The fact is that the distribution of these $\lambda_k(w)$, is exactly the 
same as the distribution of the eigenvalues of a random Hermitian $n$ by 
$n$ matrix, suitably normalized, as $n$ goes to infinity.  Does this have any 
connection with the theorems presented here?  

It may be worth pointing out explicitly that although the problems solved in 
this story range over several areas of mathematics --- 
including linear algebra, 
commutative algebra, representation theory, intersection theory, and 
combinatorics ---  none of the people involved in the recent success came to 
the problems from any of these fields.  Klyachko came from studying vector 
bundles on toric varieties, Totaro from studying filtered vector spaces 
using geometric invariant theory, Knutson, 
Agnihotri and Woodward came from symplectic geometry, Tao from 
harmonic analysis, and 
Belkale from the study of local systems on Riemann surfaces.  Klyachko's 
article [Kl1] discusses some other interesting topics, such as spaces of 
polygons and Hermitian-Einstein metrics, that are related to the eigenvalue 
problem.

\enddocument